\documentclass[preprint,12pt]{elsarticle}

\usepackage{amssymb}

\usepackage[T1]{fontenc}
\usepackage[utf8]{inputenc}
\usepackage{amsmath}
\usepackage{hyperref}

\newcommand{\im}{\mathrm{i}}

\usepackage{lipsum}
\makeatletter
\def\ps@pprintTitle{%
 \let\@oddhead\@empty
 \let\@evenhead\@empty
 \def\@oddfoot{}%
 \let\@evenfoot\@oddfoot}
\makeatother

\begin{document}

\begin{frontmatter}



\title{VEXPA: Validated EXPonential Analysis through regular sub-sampling}


\author{Matteo Briani\footnote{Research supported by the Instituut
voor Wetenschap en Technologie - IWT}, 
Annie Cuyt, Ferre Knaepkens\footnote{Research supported by the Fund 
for Scientific Research - Flanders (FWO-Vlaanderen)} 
and Wen-shin Lee}
\ead{\{annie.cuyt,ferre.knaepkens,wen-shin.lee\}@uantwerpen.be}
\address{Department of Computer Science, Universiteit Antwerpen (CMI) \\
Middelheimlaan 1, B-2020 Antwerpen, Belgium\\
\ \\
Submitted 13 Sep 2017, Revised 26 Jun 2020.}

\begin{abstract}
We present a procedure that adds a number of desirable features 
to standard exponential analysis algorithms, among which output
reliability, a divide-and-conquer approach, the automatic detection of the
exponential model order, robustness against some outliers, and the possibility to parallelize the analysis.
The key enabler for these features is the introduction of uniform sub-Nyquist
sampling through decimation of the dense signal data.
We actually make use of possible aliasing effects to recondition the
problem statement rather than that we avoid aliasing.

In Section 2 the standard exponential analysis is described, including a sensitivity analysis. 
In Section 3 the ingredients for the new approach are collected, of which good use is made in Section 4 where we essentially bring everything together in what we call VEXPA.

Some numerical examples of the new procedure illustrate in Section 5 that the additional features are indeed realized and that VEXPA is a valuable add-on to any stand-alone exponential analysis. 
While returning a lot of additional output, it maintains a favourable 
comparison to
the CRLB of the underlying method, for which we here choose a matrix
pencil method. Moreover, the output reliability of VEXPA is similar to that of atomic norm
minimization, whereas its computational complexity is far less.

\end{abstract}


\end{frontmatter}



\section{Introduction}
Many real-time experiments involve the measurement of signals which fall exponentially with time.
The task is then to determine from these measurements the number of terms $n$ and the value of all the parameters in the exponentially damped model
\begin{equation}
  \phi(t) = \sum_{i=1}^n \alpha_i \exp(\mu_i t), \qquad \alpha_i, \mu_i \in \mathbb{C}. 
  \label{intromodel}
\end{equation}
In general, parametric methods as well as 
nonparametric methods sample at a rate dictated by the Shannon-Nyquist theorem \cite{Ny:cer:28,Sh:com:49}, which states that the sampling rate needs to be at least twice the maximum bandwidth of the signal.
A coarser time grid than dictated by the theory of Nyquist and Shannon causes aliasing, mapping higher frequencies to lower ones in the analysis.
We present a parametric method that samples at a rate below the Shannon-Nyquist one, while maintaining a regular sampling scheme.
The new technique is actually exploiting aliasing, to influence the
numerical conditioning of the problem statement, rather than avoiding it. 
The latter is a useful feature as parametric methods are inherently 
more sensitive to noise. Methods that achieve much better
reliability under noise, can on the other hand be sensitive to the estimated model order or
require lots of computation time \cite{denoising,fastANM}. Another
feature of the newly proposed method is that it automatically and
concurrently provides 
a quite reliable estimate of the model order $n$.

As a consequence of the lower sampling rate it is possible to perform several independent analyses over the original set of samples, each analysis starting from a decimated dataset. 
If desired, these analyses can be carried out in parallel, thus improving
the running time of the parametric method.
The independent solutions are then passed to a cluster detection algorithm in order to add a validation step to the parametric method used, a feature that is lacking in most existing implementations. 
Thanks to the possibility to work with lower sampling rates, the validation is not at the expense of additional samples. 

The decimation of the original data adds another advantage to the method,
namely the fact that the problem size of each analysis is much smaller
since it is in size reduced by the decimation factor. 
Performing the
analysis on the different decimated sets creates a divide-and-conquer
flavour which greatly improves the overall computational complexity, 
even when not executed in parallel.   

Making use of the link between Prony-based algorithms and Pad\'e
approximation, we are able to separate the uncorrelated noise from the actual signal and avoid the computation of bogus terms in case of a low signal-to-noise ratio. 
In this way the proposed method detects the number of components $n$ automatically. 
The latter is a nice side result of working with independent decimations
of the given signal data. We emphasize that our aim is not to merely
obtain an estimate for the model order $n$,
such as can be provided by information theoretic criteria
(AIC, MDL, etc.), but to actually extract the correct sparsity $n$
from the data samples. The goal is not to fit a least complex exponential model
to the data, 
but to solve the inverse problem of deducing the correct model order $n$,
which is known to be a difficult problem.

Each decimated set of samples is now subject to an independent realization of the noise. 
While an unfiltered outlier may skew a single analysis, independent decimations indicate the presence of an outlier.
The cluster analysis makes the underlying exponential analysis algorithm more 
robust with respect to such persistent outliers, which is another desirable feature.

\section{The multi-exponential model}
Exponential analysis is an inverse problem and may therefore 
be more sensitive to
noise. Besides recalling the basic theory and its connections to some
other topics, we also discuss its susceptibility to noise. 

\subsection{Exponential analysis}
Let $\phi(t)$ be a sum of complex exponentials with $\Re(\mu_i)$,
$\Im(\mu_i)$, $|\alpha_i|$ and $\arg(\alpha_i)$ respectively denoting the
damping, frequency, amplitude and phase in each component of the signal 
\begin{equation}
  \phi(t) = \sum_{i=1}^{n} \alpha_i \exp(\mu_i t), \label{model} 
\end{equation}
where the $\mu_i$ are assumed to be mutually distinct. 
We sample the function $\phi(t)$ at the points $j\Delta$ for $j = 0, \ldots,2n -1, \ldots, N-1 $ and we set $\Omega = 1/\Delta$.  
Furthermore, we assume that the frequency content $\Im(\mu_i), i=1, \ldots, n$ in $\phi(t)$ is limited by 
\begin{equation}
 |\Im(\mu_i) / (2 \pi)| < \Omega/2, \qquad i = 1,\ldots,n. 
  \label{nyquist}
\end{equation}
The aim is to extract the model order $n$ and the parameters $\mu_1,
\ldots, \mu_n$ and $\alpha_1, \ldots, \alpha_n$ from a limited number of
samples 
of $\phi(t)$.
When the data are noisefree, the $2n$ parameters $\alpha_i$ and $\mu_i$
can be extracted from $2n$ consecutive samples \cite{dePr:ess:95}. 
In order to confirm or reveal the value of $n$ at least one more sample is required \cite{Ka.Le:ear:03}. 
In a noisy context preferably more than the minimal number of samples is provided.

In the sequel we write
\begin{align*}
  \phi_j &:= \phi(j\Delta), \qquad j=0, \ldots, N-1,  \qquad N \ge 2n, \\
  \lambda_i &:= \exp(\mu_i \Delta), \qquad i=1, \ldots, n,
\end{align*}
and for integer values $s$ and $u$, we denote by
\begin{equation}
  {_u^s}H_n :=
  \begin{pmatrix}
    \phi_s & \dots & \phi_{s+(n-1)u} \\
    \vdots & \ddots & \vdots \\
    \phi_{s+(n-1)u} & \dots & \phi_{s+(2n-2)u} \\
  \end{pmatrix},
  \qquad s \geq 0, u \geq 1,
\label{hankelmatrix}
\end{equation}
the square Hankel matrix of size $n$ constructed from the samples $\phi_j$. 
The left subscript $u$ and left superscript $s$ are respectively called the {\it undersampling} and the {\it shift} parameters. 
Whenever attached to the left of a mathematical notation in the sequel, they need to be interpreted as such.

In the standard case $u=1$ and $s=0$ or $1$. 
Note that the Hankel matrices ${_1^0}H_n$ and ${_1^s}H_n$ can be decomposed as
\begin{align*}
  {_1^0}H_n &= V_n A_n V_n^T, \qquad {_1^s}H_n = V_n \Lambda_n^s A_n V_n^T, 
\\ \\
  V_n &= \begin{pmatrix}
    1 & 1 & \cdots  & 1 \\
    \lambda_1 & \lambda_2 & \cdots & \lambda_n \\
    \vdots & \vdots & & \vdots \\
    \lambda{_1^{n-1}} & \lambda{_2^{n-1}} & \cdots & \lambda{_n^{n-1}}
  \end{pmatrix},
  \qquad
  \begin{aligned}
    A_n &= \text{diag}(\alpha_1, \ldots, \alpha_n), \\
    \Lambda_n &= \text{diag}(\lambda_1, \ldots, \lambda_n).
  \end{aligned}
\end{align*}
Then the model order $n$, the coefficients $\alpha_i$ and the parameters
$\mu_i$ are retrieved from the samples $\phi_j$ using a variant of Prony's
method \cite{Ro.Ka:esp:89,Hu.Sa:mat:90,St.Yi.ea:sta:94}. 
Prony's method consists of two stages: first the parameters $\lambda_i$ are retrieved from which the $\mu_i$ can be extracted because of \eqref{nyquist}, and then the $\alpha_i$ are computed from a linear system of equations. 
Often the $\lambda_i$ are obtained from the generalized eigenvalue problem
\cite{Hu.Sa:mat:90}
\begin{equation}
  ({_1^1}H_n) v = \lambda ({_1^0}H_n) v. 
  \label{generalized_eigenvalue}
\end{equation}
Subsequently the $\alpha_i$ are computed from the interpolation conditions
\begin{equation}
  \sum_{i=1}^n \alpha_i \exp( \mu_i j\Delta) = \phi_j, \qquad j = 0, \ldots, 2n - 1,
\ldots, N-1 
  \label{vandermondealpha}
\end{equation}
either by solving the system in the least squares sense, in the presence of noise, or by solving a subset of $n$ interpolation conditions in case of a noisefree $\phi(t)$.
Note that $\exp(\mu_i j\Delta) = \lambda_i^j$ and that the coefficient matrix of \eqref{vandermondealpha} is therefore a Vandermonde matrix.
In a noisy context the Hankel matrices in \eqref{generalized_eigenvalue} can also be extended to rectangular matrices and the generalized eigenvalue problem can be considered in a least squares sense \cite{Ch.Go:gen:06}.

Condition \eqref{nyquist} guarantees that the $\mu_i$ can be extracted from the $\lambda_i$ without ambiguity. 
However, when $|\Im(\mu_i) / (2 \pi)| \geq \Omega/2$, then 
each computed $\lambda_i$ represents an entire set of possible $\mu_i$ and 
$\Im(\mu_i)$ may be identified with a smaller frequency, an effect known as {\it aliasing}. 
How to solve the aliasing problem in that case is addressed in \cite{Cu.Le:how:18} and recalled in Section 3. 

What can be said about the number of terms $n$ in \eqref{model}, which is also called the {\it sparsity}? 
From \cite[p.~603]{He:app:74} and \cite{Ka.Le:ear:03} we know that
\begin{align*}
&\det {^s_1}H_\nu = 0 \text{ accidentally}, \qquad \nu<n, \\
&\det {^s_1}H_n \neq 0, \\ 
&\det {^s_1}H_\nu = 0, \qquad  \nu > n.
\end{align*}
While the second and third statement are clear, we briefly explain
the first one. Because of the matrix factorisation of ${^s_1}H_n$ we know
that $\det {^s_1}H_n$ is a polynomial expression in terms of the
$\alpha_1, \ldots, \alpha_n, \lambda_1, \ldots, \lambda_n$. For $\nu < n$,
this expression is nonzero in general, unless the expression $\phi(t)$ and
the sample points $j\Delta$ are such that one accidentally hits a zero of
this polynomial. A simple example makes this crystal clear.
Consider
$$\phi(t) =  2\exp(\im\pi/4\, t) - \exp(\ln(2)/2\, t) - \exp\left(
(\ln(2)/2 +\im \pi/2)t \right)$$
with $\Delta=1$. Then while $n=3$, we find with $\nu=1, 2$
that $\det {^0_1}H_1 = 0$ and $\det {^0_1}H_2 = 0$. 

A standard approach to make use of these  three statements is to compute a singular value decomposition of the Hankel matrix $^0_1H_\nu$ for increasing values of 
$\nu>n$ and apply some thresholding.
In the presence of noise and/or very similar eigenvalues, this technique
is known to be unreliable \cite{Cu.Ts.ea:fai:18}.
The method proposed in Section 4 allows to automatically detect $n$ while processing the samples $\phi_j$ without having to resort to a separate singular value decomposition of $^0_1H_\nu$.

\subsection{The Pad\'e and Froissart connections}
\begin{figure}
  \centering
  \includegraphics[scale=0.20]{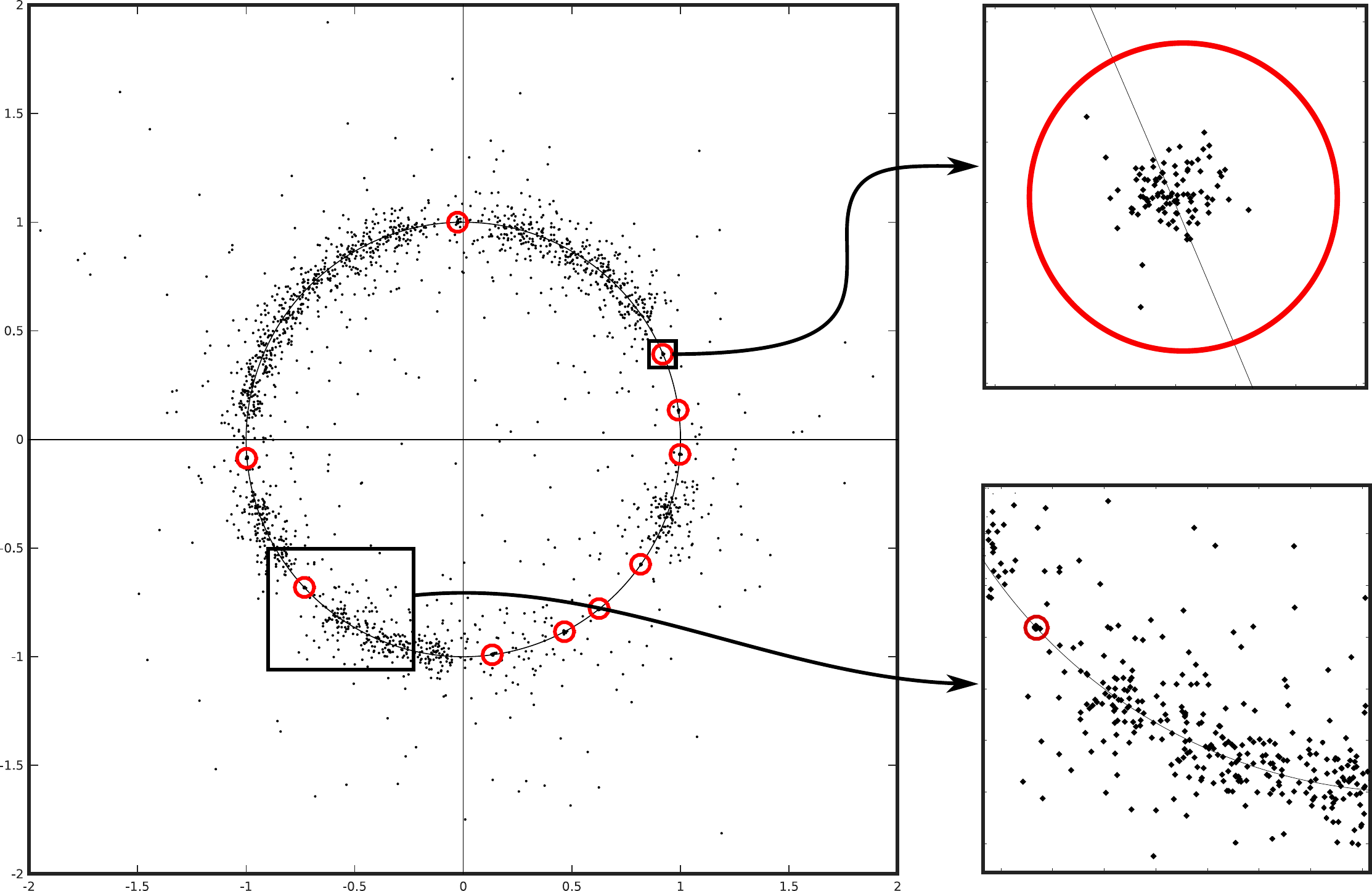}
  \caption{\it Typical analysis result ($n=10, \nu=30$) of a test signal $\phi(t)$ after several perturbations $\epsilon(z)$: the true $\lambda_i$ are drawn as red circles.}
  \label{noisebehaviour}
\end{figure}

There is an interesting but somewhat unknown connection between Pad\'e
approximation, Froissart doublets and the Prony problem, which we briefly
recall from \cite{We.Mc:pro:63,Ba.My.ea:pad:89}. 
Consider  the function $f(z)$ defined by
\begin{equation*}
  f(z) = \sum_{j=0}^\infty \phi_j z^j.
\end{equation*}
For $\phi_j=\phi(j\Delta)$ with $\phi(t)$ given by \eqref{model}, we can write
\begin{equation}
    f(z) = \sum_{i=1}^n {\alpha_i \over 1-\lambda_i z}.
\label{laplace}
\end{equation}
The partial fraction decomposition \eqref{laplace} is related to both the
Laplace transform and the Z-transform of \eqref{model} as described in
\cite{We.Mc:pro:63,Ba.My.ea:pad:89}.
It is a rational function of degree $n-1$ in the numerator and degree $n$ in the denominator with poles $1/\lambda_i$.
Now let us perturb $f(z)$ with white circular Gaussian noise to obtain 
\begin{equation*}
  f(z) + \epsilon(z) = \sum_{j=0}^{\infty} (\phi_j + \epsilon_j) z^j.
\end{equation*}
The theorem of Nuttall-Pommerenke states that if $f(z) + \epsilon(z)$ is analytic throughout the complex plane, except for a countable number of poles \cite{Nu:con:70} and essential singularities \cite{Po:pad:73}, then its sequence of paradiagonal Pad\'e approximants $\{ r_{\nu-1,\nu}(z) \}_{\nu \in \mathbb{N}}$ of degree $\nu - 1$ over $\nu$ converges to $f(z) + \epsilon(z)$ in measure on compact sets.
This means that for sufficiently large $\nu$ the measure of the set where the convergence is disrupted, so where $|f(z) + \epsilon(z) - r_{\nu -1, \nu} (z)| \geq \tau$ for some given threshold $\tau$, tends to zero as $\nu$ tends to infinity.
Pointwise convergence is disrupted by $\nu-n$ unwanted pole-zero
combinations of the Pad\'e approximants that are added to the $n$ true
poles and $n-1$ true zeros of $f(z)$ \cite{Ga:eff:72}, with the pole and
zero in the undesirable pair almost cancelling each other locally.  
These pole-zero combinations are also referred to as Froissart doublets.  
In practice, these Froissart doublets offer a way to separate the noise
$\epsilon(z)$ from the underlying $f(z)$ \cite{Be:pad:96}.
Because of the Pad\'e convergence theorem, the true (physical) poles can be identified as stable poles in successive $r_{\nu-1,\nu}(z)$, while the spurious (noisy) poles are distinguished by their instability.
When increasing $\nu$ we compute a larger set of poles, of which the noisy ones are moving around in the neighbourhood of the complex unit circle \cite{Gi.Pi:pad:97,Gi.Pi:pad:99} with every different realization of the noise $\epsilon(z)$. 
The latter is illustrated in Figure \ref{noisebehaviour} where we show the results of the analysis of a test signal perturbed by a large number of independent noise realizations: the true $\lambda_i$ are forming clusters while the ones related to noise are scattered around \cite{BARONE2005224,perotti2014identification}.
In addition, around each $\lambda_i$-cluster one empirically finds an almost Froissart doublet-free zone.

This characteristic of the true poles is a key point on which our method is based: after the computation of $\nu>n$ generalized eigenvalues $\lambda_i$, we discard the unstable ones and focus on the stable ones.
We now describe in more detail the precise influence of noise in the data $\phi_j$ on the $\lambda_i$.

\subsection{Sensitivity to noise}
The exponential analysis of $\phi(t)$, being an inverse problem, is known
to be sensitive to noise. Here we briefly recall what is known and in the
next section we explain how the new method is able to deal with certain outliers on the one hand and normally distributed noise on the other. 

In \cite{Go.Mi.ea:sta:99} the authors explain that the roundoff errors in the computation of the generalized eigenvalues are amplified by mainly three sources:
\begin{itemize}
  \item the scaling of the problem (the $\lambda_i$ should lie as closely as possible to the complex unit circle),
  \item the size of the $|\alpha_i|$ relative to the noise ($\lambda_i$ with smaller amplitude are more challenging to retrieve),
  \item the relative position of the $\lambda_i$ with respect to each
other (clustered $\lambda_i$ are hard to separate and retrieve
individually).
\end{itemize}
The first problem is addressed in \cite{Go.Mi.ea:sta:99} by means of a
diagonal preconditioning matrix, and in \cite{Be.Go.ea:num:07} by using a suitably chosen invertible upper triangular matrix.
The second problem can be tackled with the use of linear time invariant filters which preserve model \eqref{model}.
A solution for the third problem is proposed in \cite{Cu.Le:how:18} and accomplishes a redistribution of the $\lambda_i$. 
Our new method is based on this approach.
We now briefly recall the basics of the analysis in \cite{Be.Go.ea:num:07}
to understand the effect of noise and how this is related to the method
presented in \cite{Cu.Le:how:18}. 

\begin{figure}
  \centering
\includegraphics[height=6truecm, width=6truecm]{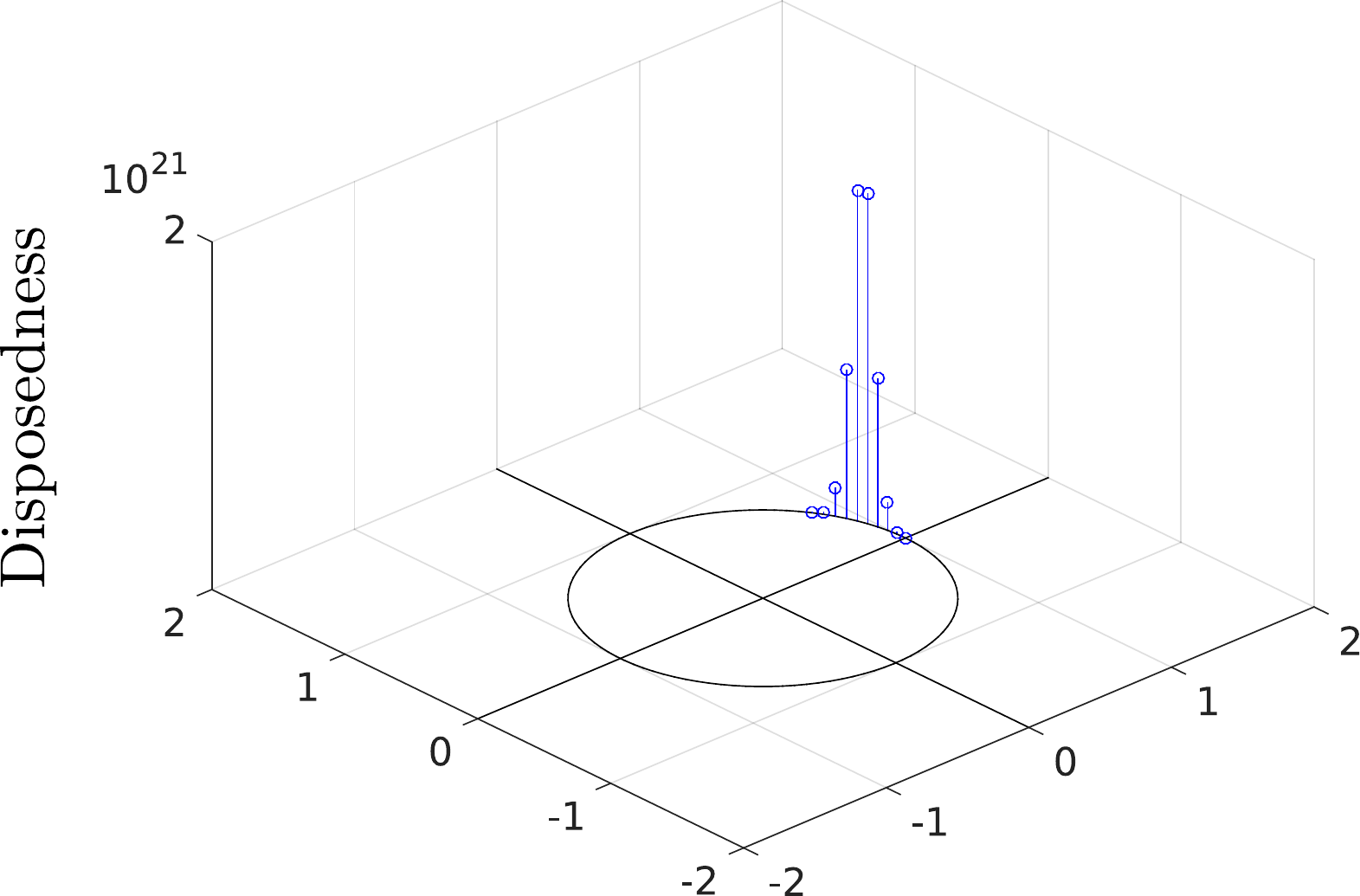}\quad
\includegraphics[height=6truecm, width=6truecm]{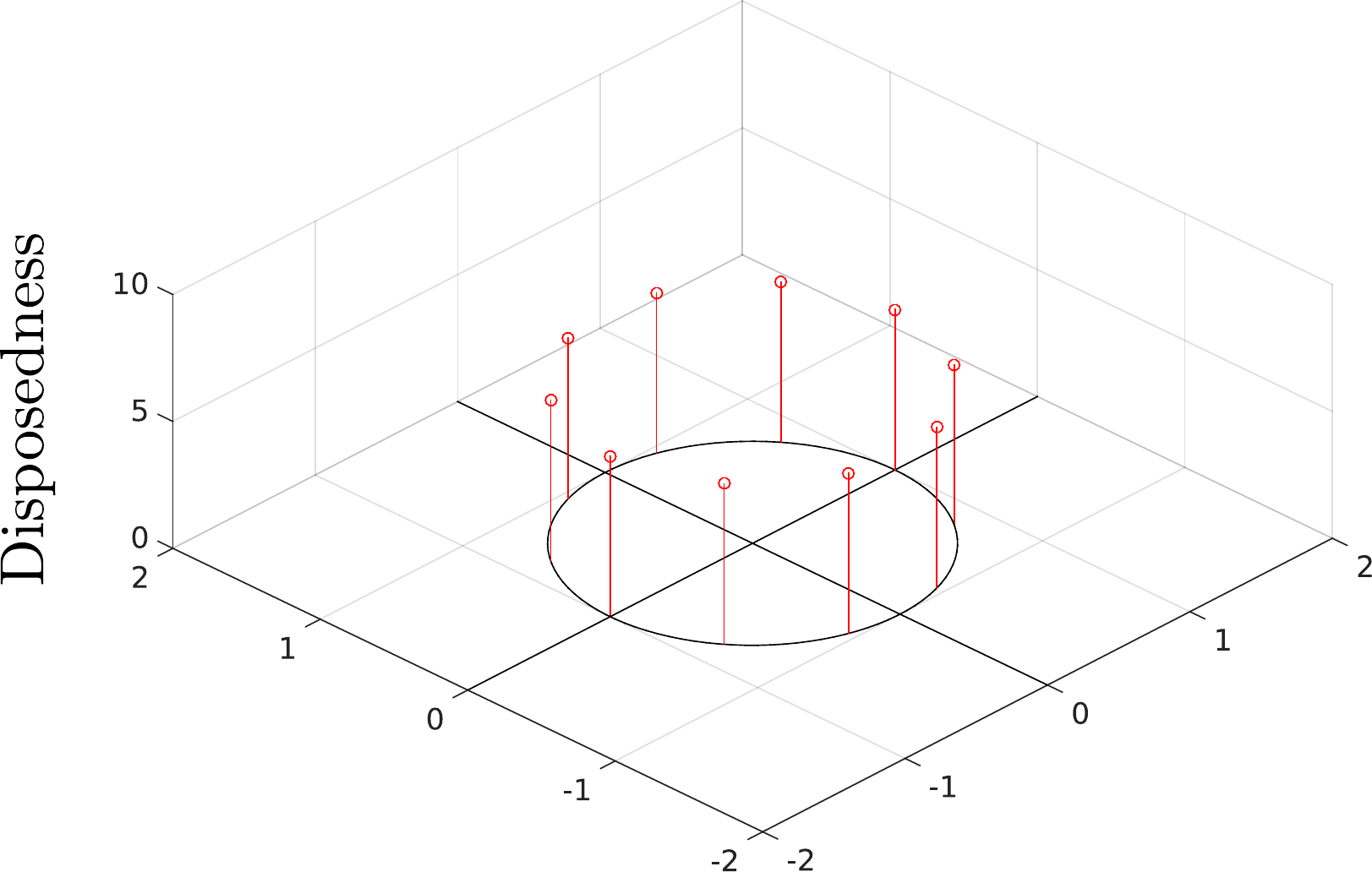} 
\caption{\it Ill-disposed $_1\lambda_i=\exp(\im 2\pi(i-1)/100)$ at the left and well-disposed $_{10}\lambda_i=\exp(\im 2 \pi (i-1)/10)$ at the right, $i=1, \ldots, 10$.}
\label{aliasing_toy_example}
\end{figure}

Let $(\epsilon_0, \ldots, \epsilon_{2n-1}, \ldots, \epsilon_{N-1})$ again denote the noise vector added to the samples $(\phi_0, \ldots, \phi_{2n-1},$ $\ldots, \phi_{N-1})$. 
We rewrite the noise terms $\epsilon_j$ as $\epsilon_j=\epsilon e_j$ where the square Hankel matrices $^0_1E_n$ and $^1_1E_n$ of size $n$, filled as in \eqref{hankelmatrix} but now with the $e_j$ instead of the $\phi_j$, satisfy
\begin{equation*}
  ||{^0_1}E_n||_2 \le 1, \qquad ||{^1_1}E_n||_2 \le 1.
\end{equation*}
Let $L_i(\lambda)$ denote the Lagrange basis polynomial of degree $n-1$ with
roots $\lambda_1, \ldots$, $\lambda_{i-1}$, $\lambda_{i+1}$,$ \ldots,
\lambda_n$ and $L_i(\lambda_i)=1$, so
$$L_i(\lambda) = \frac{\prod_{k=1, k\not= i}^n
(\lambda-\lambda_k)}{\prod_{k=1, k\not= i}^n
(\lambda_i-\lambda_k)}.$$ 
The coefficients of the polynomial
$L_i(\lambda)$ make up the vector $\ell_i$ of size $n$.  
Then the {\it disposedness} $\rho_i$ of the  generalized eigenvalue 
$\lambda_i(\phi_0+\epsilon e_0, \ldots, \phi_{N-1}+\epsilon e_{N-1})$,
as a function of the given $\phi_j$ and the noise terms
$\epsilon_j=\epsilon e_j$, is defined by
\begin{equation*}
  \rho_i := \left |{d \lambda_i \over d \epsilon}(0) \right |
\end{equation*}
and satisfies
\begin{equation}
\rho_i \leq {|\lambda_i| + 1 \over |\alpha_i|} \; ||\ell_i||_2^2 \; ( ||{_1^1}H_n||_2  + ||{_1^0}H_n||_2 ). 
  \label{first_approx}
\end{equation}
A generalized eigenvalue $\lambda_i$ is ill-disposed when $\rho_i$ is large. 
Larger $\rho_i$ imply higher susceptibility to noise. 
Besides the Froissart phenomenon described earlier, the disposedness
$\rho_i$ of the generalized eigenvalues, or rather its computable upper
bound given in \eqref{first_approx}, 
is another tool to use when inspecting the $\lambda_i$.  
In Figure \ref{aliasing_toy_example} we illustrate the relationship between the $\rho_i$ and the
relative position of the $\lambda_i$ with respect to each other: we plot
the right hand side of \eqref{first_approx} for a toy problem where  
we choose $\Omega=100, n=10, \alpha_i=1, \mu_i= \im 2 \pi (i-1)$.  
At the left the upper bounds for the values $\rho_i$ are plotted at the
locations of the generalized eigenvalues $\lambda_i = \exp(\mu_i\Delta) =
\exp(\mu_i/\Omega), i=1, \ldots, 10$. Now let us
change the undersampling parameter $u$ in $^s_uH_n$ in
\eqref{hankelmatrix} and \eqref{first_approx} from $u=1$ to $u=10$, which
is equivalent to replacing $\Delta$ by $u\Delta$ or replacing $\Omega$ by
$\Omega/u$. We
recompute the generalized eigenvalues $_u\lambda_i = \exp(\mu_i
(10\Delta)) = \exp(10\mu_i/\Omega)$ and the
disposedness, which we now denote by $_u\rho_i$. The result, which is
shown at the right, changes dramatically, from $O(10^{21})$ to almost
$O(10^1)$. 

Another important tool for inspecting the $\lambda_i$ is the Cram\`er-Rao lower bound (CRLB) \cite{Kay:1993:FSS:151045, Yao1995CramerRaoLB}.
For any given unbiased estimator of the parameters in \eqref{model} and a specific amount and type of noise, the CRLB returns the minimal variance that the estimator suffers.
In our case, the estimator is any implementation of Prony's method and the type of noise is white circular Gaussian noise.
The CRLB depends on the number of samples $N$, the variance and type of
noise and the set of parameters $|\alpha_i|, \arg(\alpha_i), \Re(\mu_i)$ and
$\Im(\mu_i)$. 
The bound is often used to compare the variance of a specific estimator to this theoretical lower bound.
The closer an estimator is to the CRLB, the more efficient it is said to be.

\begin{figure}
\centering
\includegraphics[scale=0.20]{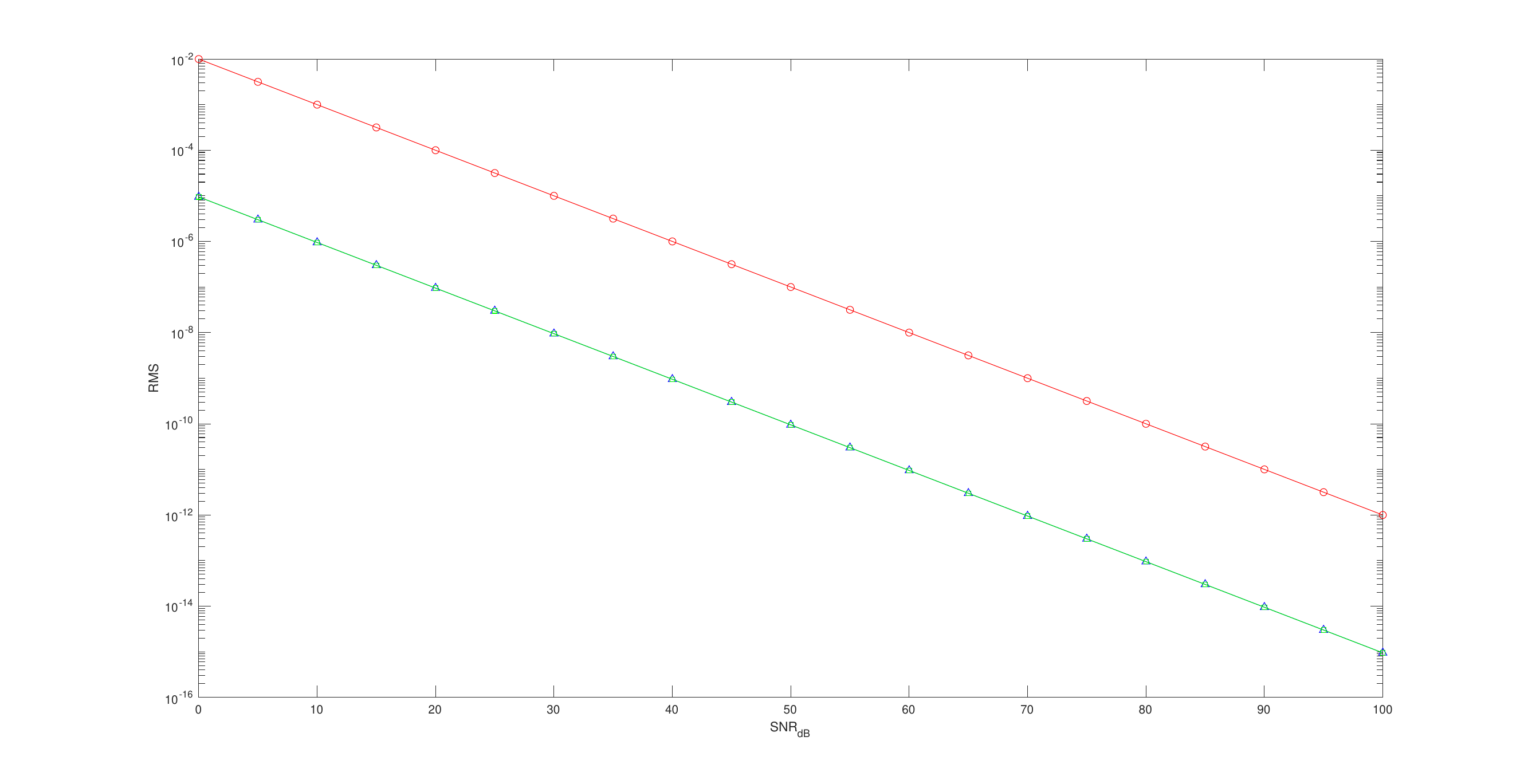}
\caption{\it Root mean square of the CRLB vector of the $\Im(\mu_i), i = 1,\ldots,10$, respectively for $\Omega=100, N=200$ (blue), $\Omega=10, N=200$ (green), $\Omega=10, N=20$ (red).}
  \label{CRLBexample}
\end{figure}

We consider the practical computation of the CRLB provided in \cite{Yao1995CramerRaoLB} and illustrate the relationship between the CRLB and the disposedness $\rho_i$ of $\lambda_i, i=1, \ldots, n$.
Take the same toy example and add white circular Gaussian noise of varying signal to noise ratio (SNR). 
In Figure \ref{CRLBexample} we graph the root mean square of the vector of CRLB's for the
parameters $\Im(\mu_i), i=1, \ldots, 10$, and this for decreasing SNR in three different situations:
\begin{itemize}
  \item $\Delta=1/\Omega, N=200$ samples $\phi_j$ (blue triangles),
  \item $\Delta=10/\Omega, N=200$ samples $\phi_j$ (green squares),
  \item $\Delta=10/\Omega, N=20$ samples $\phi_j$ (red circles).
\end{itemize}
Note that multiplying $\Delta$ by $u=10$ while maintaining $N=200$ implies that the signal is sampled over a larger time interval, while multiplying $\Delta$ by $u=10$ and dividing $N$ by $u=10$ does not enlarge the observation window.
So in the first and second case the number of samples is equal while in the
first and third case the observation window is equal. Our aim is to
get the best of both worlds: while decimation of the signal samples takes
you from the CRLB in blue (triangles) to the CRLB in red
(circles), we want to recombine separate decimations in order to return
from the latter to the former while profiting from some
additional features on the way. Decimation significantly diminishes the size
of the generalized eigenvalue problems, improves the numerical conditioning, 
and will automatically 
return a reliable estimate for the model order $n$.
How this can be done is described in Section 4 and illustrated in Figure
\ref{fig:CRLB}. 

\section{Recovering from aliasing after decimation}
So we know that choosing $u > 1$ may positively impact the disposedness of the $\lambda_i$, without negatively impacting the CRLB if the total number of samples can approximately be maintained. 
Since introducing $u$ impacts $\Delta$ or $\Omega$, aliasing may occur when \eqref{nyquist} is violated. 
We now explain how to deal with this effect: the goal is to enjoy the positive influence of a larger $u$ without suffering the aliasing effect introduced by it. 
\subsection{Decimation}
\begin{figure}
  \centering
\includegraphics[scale=0.35]{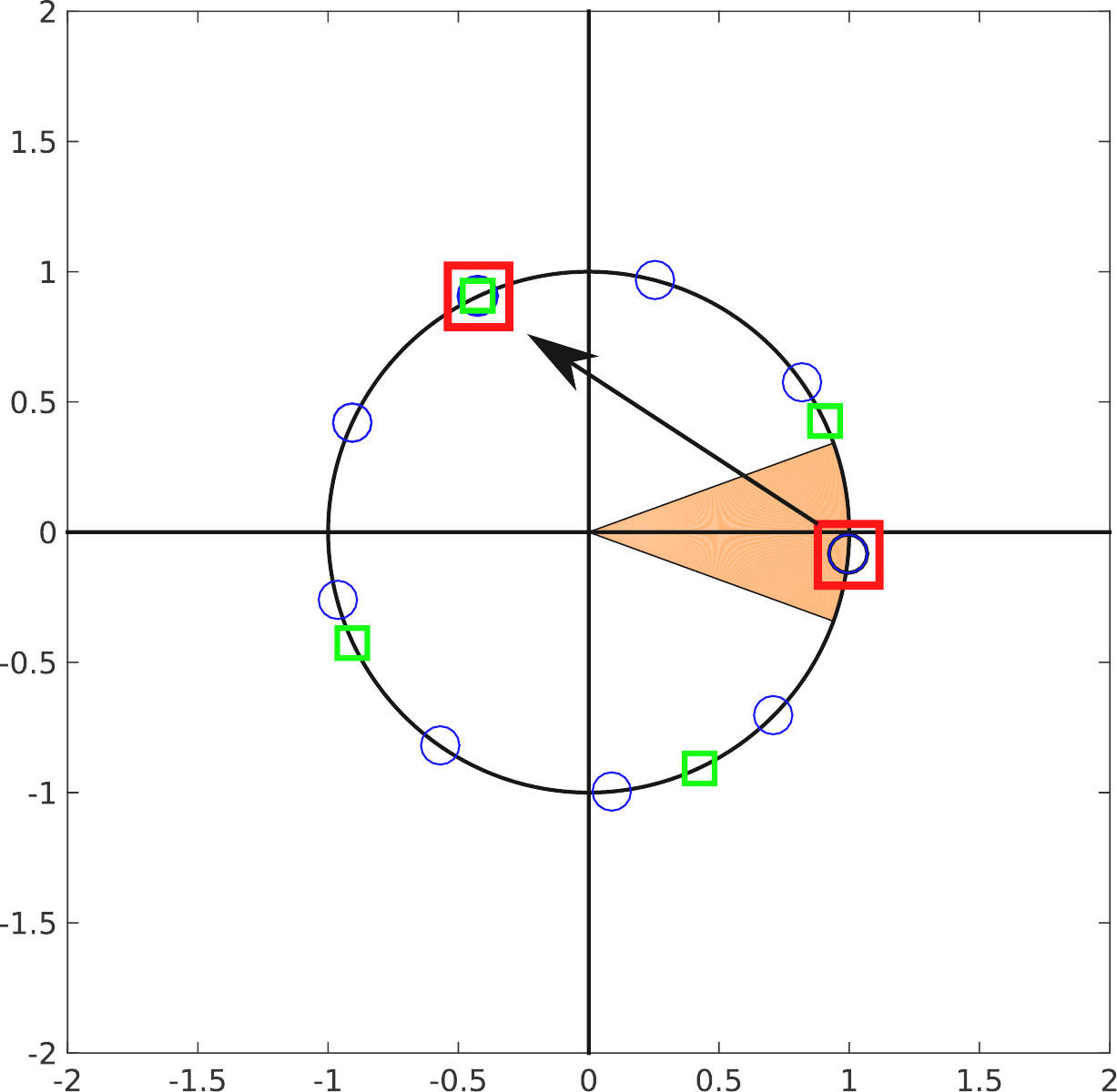}
\caption{\it Intersection of $U_i$ (blue circles, $u=9$) and $S_i$ (green squares, $s=4$), relocating the aliased ${_u\lambda}_i$ (red square).}
 \label{roots_intersection}
\end{figure}

Instead of using the consecutive set of samples $\phi_j, j=0,\ldots,2n-1, 
\ldots N-1$, we consider the decimated set $\phi_{uj}$ which is obtained by considering one sample every $u$ samples, thus sampling $\phi(t)$ at $j(u
\Delta)$. 
The generalized eigenvalue problem 
\begin{equation*}
  ({^u_u}H_n )v  = \lambda ({^{0}_u}H_n) v,
\end{equation*}
leads to a new set of generalized eigenvalues
\begin{equation*}
  {_u\lambda}_i := \exp(\mu_i u\Delta) = \lambda_i^u, \qquad i = 1, \ldots, n.
\end{equation*}
From ${_u\lambda}_i$ we cannot directly retrieve $\lambda_i$, due to the disruption of \eqref{nyquist}.
We are left with a set of possible values for $\lambda_i$ given by
\begin{equation*}
  U_i := \left\{ \exp\left( \mu_i\Delta + {2 \pi \im \over u} \ell
\right), \ell=0, \ldots, u-1 \right\}.
\end{equation*}
Despite this, we can already compute the coefficients $\alpha_i$ by solving the linear system
\begin{equation}
  \phi_{uj} = \sum_{i=1}^{n} \alpha_i ({_u\lambda}_i)^j, \qquad j=0, \ldots, 2n-1, \ldots 
  \label{original_VDM}
\end{equation}
Now we consider a shifted set of samples $\phi_{s+uj}$ consisting of at least $n$ samples, for instance at $j=k, \ldots, k+n-1, 0 \le k \le n$, and we choose $s$ coprime with $u$.
Since
\begin{equation}
  \phi_{s+uj} = \sum_{i=1}^n (\alpha_i \lambda_i^s)({_u\lambda}_i)^j, \qquad
  j=k, \ldots, k+n-1, 
  \label{undersampling_core}
\end{equation}
we denote the coefficient of $({_u\lambda}_i)^j$ in the shifted sample $\phi_{s+uj}$ by
\begin{equation*}
^s\alpha_i:= \alpha_i\lambda_i^s, \qquad i=1, \ldots, n.
\end{equation*}
We can solve the interpolation conditions \eqref{undersampling_core} for the second set of coefficients ${^s}{\alpha}_i$. 
Note that the linear systems \eqref{undersampling_core} and
\eqref{original_VDM} have the same Vandermonde structured coefficient
matrix, except for the size. This precisely connects the two coefficients
$\alpha_i$ and ${^s}\alpha_i$, and consequently ${^s}\lambda_i$, 
to the ${_u}\lambda_i$.
From $\alpha_i$ and ${^s}{\alpha}_i$ we obtain 
\begin{equation*}
  {^s}{\alpha}_i/{\alpha_i} = \lambda_i^s,
\end{equation*}
which we can denote by ${^s\lambda}_i$.  
Due to the same possible disruption of condition \eqref{nyquist}, ${^s\lambda}_i$ also stands for a set of possible values for $\lambda_i$, namely
\begin{equation*}
  S_i := \left\{ \exp\left( \mu_i\Delta + {2 \pi \im \over s} \ell
\right), \ell=0, \ldots, s-1 \right\}.
\end{equation*}
Both sets $U_i$ and $S_i$ contain the solution $\lambda_i$.
Since $u$ and $s$ are coprime they share one and only one element which is the non-aliased $\lambda_i$ \cite{Cu.Le:how:18}. 
In Figure \ref{roots_intersection} we graphically sketch what happens. 
There $u=9$, the elements in $U_i$ are shown using blue circles, $s=4$, the elements in $S_i$ are shown using green squares and the arrow points to the unique non-aliased $\lambda_i$ in their intersection. 
The orange portion is the region where the aliased ${_u\lambda}_i$ lies
(red square), from which we have to recover the correct $\lambda_i$. The
aliasing is the consequence of the decimation of the collected samples by a factor $u$.

While in theory $u$ and $s$ may be chosen arbitrarily large, this is not 
the case in practice, since noise can make it hard to point at the one 
correct value when a large number of points lie closer together in $U_i$ and 
$S_i$. Usually the chosen value of $u$ is larger than that of $s$. 
Numerical experiments indicate that a smaller $s$ is more important than a 
smaller $u$: the ${_u\lambda}_i$ values are usually less affected by noise 
than the ${^s\lambda}_i$ which are obtained as solution of Vandermonde 
structured linear systems.

\subsection{Recovery}
While we know theoretically that $U_i$ and $S_i$ have only one element in their intersection, we still need to find a way to compute this element in practice. 
In \cite{Cu.Le:how:18} the following two options are presented.  
Here we develop a more robust third approach.

An obvious approach is to compute all distances between elements of $U_i$ and elements of $S_i$ and select the pair that lies closest. 
This simple approach does not deliver satisfactory results though, because of noise issues. 
For increasing noise levels, the sets $S_i$ may be too perturbed, thus leading to a wrong match of the candidate values for $\lambda_i$.

A less obvious approach is to use the Euclidean algorithm and compute two integers $w$ and $r$ satisfying $wu+rs=1$ for the coprime $u$ and $s$. 
Then $\lambda_i$ can be retrieved as
  \begin{equation*}
    ({_u\lambda}_i)^w \; ({^s\lambda}_i)^r = \exp((wu+rs)\mu_i\Delta) =
    \lambda_i.
  \end{equation*}
The downside of this method is that if $w$ and $r$ are not small, any noise present in ${_u\lambda}_i$ and ${^s\lambda}_i$ is amplified.

We propose to solve a small number of additional systems of the form \eqref{undersampling_core}, in order to stabilize the location of the elements in $S_i$ before building the distance matrix.
We continue the use of shifted samples:
  \begin{equation}
    \phi_{ms + uj} = \sum_{i = 1}^n ({^{ms}\alpha}_i) \; ({_u\lambda}_i)^j, \qquad m = 0,\ldots, M -1. \label{pmsystems}
  \end{equation}
From each shift we compute the coefficients ${^{ms}\alpha}_i$ and we set up the sequence of values
  \begin{equation*}
    \alpha_i, {^s\alpha}_i, \ldots, {^{ms}\alpha}_i, \ldots, {^{(M-1)s}\alpha}_i,
  \end{equation*}
satisfying
  \begin{equation}
  {^{ms}\alpha}_i = \alpha_i ({^s\lambda}_i)^m=\alpha_i \exp(\mu_i m(s\Delta)), \qquad m=0, \ldots, M-1, 
    \label{msalpha}
  \end{equation}
where ${^0\alpha}_i = \alpha_i$.
So for fixed $i$ the values ${^{ms}\alpha}_i$ follow the exponential model
\eqref{msalpha} consisting of only one term. We can therefore use a
Prony-like method to extract ${^s\lambda}_i$ in \eqref{msalpha} from the values ${^{ms}\alpha}_i$, just as described in the previous section on basic exponential analysis.
This approach stabilizes the location of ${^s\lambda}_i= \lambda_i^s$ by the use of extra estimates. 

At this point we want to point out and stress, that the whole procedure of decimation and recovery can be used on top of any Prony-like method. 
Retrieving ${_u\lambda}_i, {^s\lambda}_i$ or ${^{ms}\alpha}_i$ for chosen $u$ and $s$ does not require a specific parametric method. 
In fact, the current procedure offers a way to parallelize existing Prony-like methods, as the decimated signals can be treated independently of each other.
In the next section we explain how the combination of the decimated
results adds, as one of the features, a validation step to the method, which is mostly lacking in existing Prony-like algorithms.

\subsection{Frequency collision}
A problem that may occur when decimation causes aliasing, is the possible collision of frequencies. For instance, two distinct eigenvalues $\lambda_1$ and $\lambda_2$ may be aliased to the same eigenvalue ${_u\lambda}_1={_u\lambda}_2$.
However unlikely, we want to discuss how to deal with this situation. We
explain the remedy on an example. A fully detailed mathematical analysis
of all the eventualities that can occur as a consequence of the
decimation, is presented in \cite{Cu.Le:how:18}.

Let $\phi(t)$ be specified by $n = 2$, $\alpha_1 =\alpha_2 = 1$, $\mu_1 = 2 \pi \im 13$, $\mu_2 = 2 \pi \im 33$. 
We set $\Omega = 100$ and consider one sample $\phi_j=\phi(j/\Omega)$ every ten samples ($u=10$) thus changing $\Omega$ to be $10$.
Due to aliasing, $\lambda_1$ and $\lambda_2$ are mapped to another location in the complex plane.
In particular, we have
  \begin{equation*}
    {_u\lambda}_1 = {_u\lambda}_2 = \exp\left({2 \pi \im 3 \over 10} \right)
  \end{equation*}
because
  \begin{equation*}
\exp\left({2 \pi \im 33 \over 10}\right) = \exp\left({2 \pi \im 13 \over
10}\right) = \exp\left({2 \pi \im 3 \over 10}\right).
  \end{equation*}
So in the decimation step \eqref{original_VDM} Prony's method retrieves a single frequency with associated coefficient $\alpha_1+\alpha_2$. 

It is however still possible to retrieve the original values $\lambda_1$ and $\lambda_2$ in the recovery step.  
As explained, the generalized eigenvalue ${_u\lambda}_1={_u\lambda}_2$ stands for a set of values $U_1=U_2$ that now contains both the correct $\lambda_1$ and $\lambda_2$. 
We choose $s$ coprime with $u$ and compute the values ${^{ms}\alpha}_1$ (remember that the computed ${^0\alpha}_1=2$ now equals the sum of the true coefficients).
Since $s$ is coprime with $u$, no frequency collision occurs in ${^{ms}\alpha}_1$ which is following the model
\begin{equation}
  {^{ms}\alpha}_1 = \alpha_1 \exp(\mu_1 m s\Delta) + \alpha_2 \exp(\mu_2 m s \Delta), \qquad m=0, \ldots, M-1. 
  \label{msalpha_coll}
\end{equation}
So in the analysis of \eqref{msalpha_coll} Prony's method reveals two
contributions ${^s\lambda}_1$ and ${^s\lambda}_2$ which bring forth the
sets $S_1$ and $S_2$, respectively containing $\lambda_1$ and $\lambda_2$.
The intersections $U_1\cap S_1$ and $U_2 \cap S_2=U_1\cap S_2$ reveal the original $\lambda_1$ and $\lambda_2$. 

\begin{figure}
\centering
\includegraphics[scale=0.29]{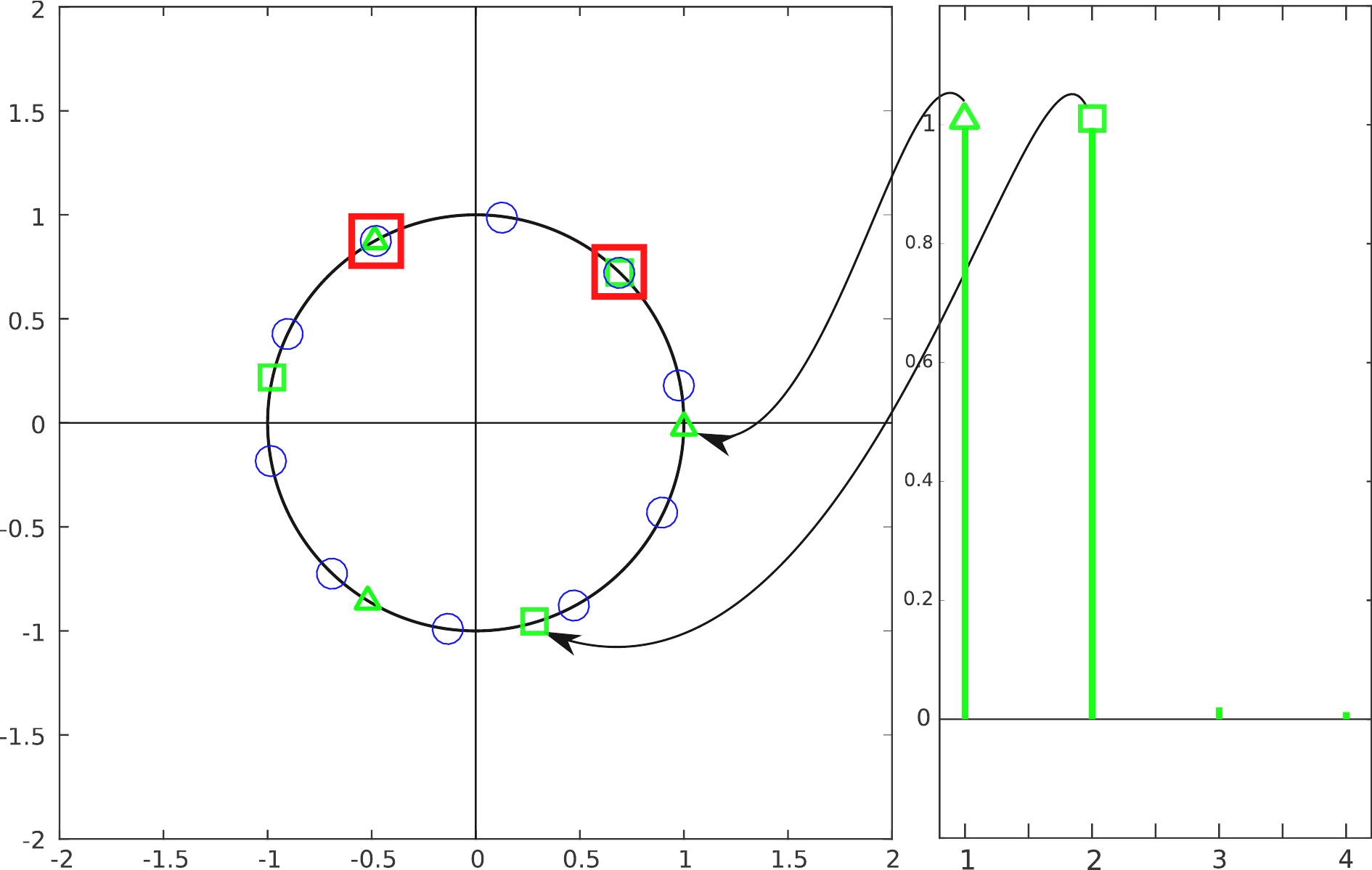}
\caption{\it The $|\alpha_1|$ of $\;{^s\lambda_1}$ (green triangle) and $|\alpha_2|$ of $\;{^s\lambda}_2$ (green square) at the right ($s=3, M=8$), identifying $U_1\cap S_1=\{\lambda_1\}$ and $U_2\cap S_2=U_1\cap S_2=\{\lambda_2\}$ from $u=10$ at the left (red squares).}
  \label{fig:collision}
\end{figure}

Of course the above can also be applied to the more general case of several collisions in a signal $\phi(t)$ containing more terms. 
The key element is that the value $M$ in \eqref{msalpha_coll} is chosen large enough to allow the identification of all the collided eigenvalues. 
In particular, $M$ should be at least twice the number of collided eigenvalues. 
Since this number is unknown, the standard procedure is to take $M$ even and fit the ${^{ms}\alpha}_i$ with a model of size $M/2$. 
If less than $M/2$ frequencies have collided, then some of the terms in the expression for ${^{ms}\alpha}_i$ model the noise and can easily be discarded, as explained in Section 2.
We show a typical situation in Figure \ref{fig:collision}, which applies to the $n=2$
example above: the set $U_1=U_2$ is depicted using blue circles ($u=10$)
and the sets $S_1$ and $S_2$ using green triangles and squares respectively ($s=3$).
We choose $M=8$. The intersections $U_1\cap S_1$ and $U_2\cap S_2$ are indicated using red squares.

\section{Validated exponential analysis}
A quite robust Prony-like implementation, which approaches the theoretical
CRLB (depicted using blue triangles in Figure \ref{CRLBexample}), is for instance found in
\cite{Hu.Sa:mat:90,YILMAZER2008561}. 
In the sequel we refer to this method as MP, from Matrix Pencil.
Our aim now is to maintain as good as the same accuracy, but add the following features to the implementation by making a detour via decimation: 
\begin{itemize}
  \item validation of the output,
  \item automatic estimation of the model order $n$,
  \item robustness against some outliers,
  \item parallelism in the algorithm.
\end{itemize}
In other words, while the sub-sampling of a signal usually leads to cruder
estimates of the already aliased frequencies (upper CRLB curve in Figure
\ref{CRLBexample}), the method explained below still achieves the desired
CRLB curve (middle curve in Figure \ref{CRLBexample}), while adding a number of desirable features that become available through the technique described in Section 3.

Given a fixed undersampling parameter $u$, we can consider $u$ decimated
sample sets $\Phi_k, k=0, \ldots, u-1$, starting respectively at $0,
\Delta, \ldots, (u-1)\Delta$. 
The first set contains $\lfloor N/u\rfloor$ samples and all subsequent sets contain either the same number of samples or one less:
\begin{equation*}
  \Phi_k:= \{\phi_{uj+k}: j=0, \ldots, \min(\lfloor N/u \rfloor, \lfloor (N-k)/u \rfloor)-1 \}, \quad k=0, \ldots,  u-1.
\end{equation*}
From each decimated set $\Phi_k$ we extract ${_u\lambda}_i, {^s\lambda}_i,
i =1, \ldots, n$ which should carry a second index $k$ now to indicate
from which decimation $\Phi_k$ the values were obtained.  
The same holds for the coefficients $\alpha_i$. 
For the sequel we therefore introduce the notations $({_u\lambda}_{i,k}), ({^s\lambda}_{i,k}), {^{ms}\alpha}_{i,k}$ with obvious meanings.
We also introduce
\begin{equation*}
{_uL} := \cup_{i=1, k=0}^{n,u-1} \{ {_u\lambda}_{i,k} \}, \qquad
{^sL} := \cup_{i=1, k=0}^{n,u-1} \{ {^s\lambda}_{i,k} \}.
\end{equation*}
We remark that the index $i$ still runs from 1 to $n$ even if the undersampling has caused collisions.
Then some ${_u\lambda}_{i,k}$ are merely duplicated.  

Each dataset $\Phi_k$ is now a decimation of the set of samples $\{\phi_0,$ $\phi_1,$ $\ldots,$ $\phi_{N-1}\}$. 
From this section on, each sample $\phi_j$ is always perturbed by noise,
but we choose to abuse the notation $\phi_j$ instead of
$\phi_j+\epsilon_j$ in order to not overload the presentation.
Each set $\Phi_k$ is subject to an independent realization of the noise because the latter affects each decimated signal in a different and independent way. 
Thanks to the connection with the theory of Pad\'e approximation and
Froissart doublets, we know that the ${_u\lambda}_{i,k}$ and
${^s\lambda}_{i,k}$ form clusters in the sets ${_uL}$ and ${^sL}$
respectively, around the  true ${_u\lambda}_i=\lambda_i^u$ and
${^s\lambda}_i=\lambda_i^s$ with $i=1, \ldots,n$. 
Any generalized eigenvalues retrieved from overestimating the model order $n$ by $\nu>n$, model the noise and are found scattered around the complex unit disk, as explained in Section 2. 
To detect the clusters in ${_uL}$ and ${^sL}$ we propose to use the density based cluster algorithm DBSCAN \cite{Ester96adensity-based}.

DBSCAN requires two additional parameters: the density $\delta$ of the clusters and the minimum number $m_\delta$ of required cluster elements. 
These parameters are chosen in terms of the noise in the signal.  
Larger values of $\delta$ allow the detection of wider clusters, which is useful in case of a higher noise level.
Smaller values of $\delta$ allow to detect denser clusters, which appear in case of very stable estimates ${_u\lambda}_{i,k}$ or low levels of noise.
A value for $m_\delta$ smaller than $u$ allows to discard bogus estimates
appearing as a consequence of, for instance, an outlier in the data.
When $m_\delta$ is set equal to $u$, each ${_u\lambda}_i$ needs to be confirmed by all the decimated analyses.
Remember that, through the coefficient matrix shared between \eqref{original_VDM} and \eqref{undersampling_core}, each element from ${^sL}$ is connected to an element in ${_uL}$. 
So any cluster detected in ${^sL}$ is tied to a set of elements from ${_uL}$ of the same size.  
We also point out that the introduction of decimation reduces the
complexity of the numerical algorithm and parallelizes the 
exponential analysis. Instead of solving a single large 
structured generalized eigenvalue problem, one is facing $u$ much smaller 
structured generalized eigenvalue problems, which makes a big difference
even when solved sequentially.
Each $\Phi_k$ is analyzed independently and the computation of the ${_u\lambda}_{i,k}$ and ${^s\lambda}_{i,k}$ does not need data from other decimations. 
All the results are collected after the individual runs and then passed to the cluster analysis. 

Essentially three different DBSCAN scenario's can occur, as sketched in
Figure \ref{cluster_explanation}: at the left we find the result of running DBSCAN on the set ${_uL}$ and at the right the result on the set ${^sL}$.

\subsection{Standard scenario}
A cluster ${_uC}_1$ is detected in the set ${_uL}$ and its center of gravity can serve as an estimate of one of the $\lambda_i^u$. 
The elements ${^s\lambda}_{i,k}$ tied to the generalized eigenvalues ${_u\lambda}_{i,k} \in {_uC}_1$ also form a cluster, which we denote by ${^sC}_1$. 
Its center of gravity then returns an estimate of $\lambda_i^s$. 
From both centers of gravity a reliable estimate of $\lambda_i$ can be
extracted as described in Section 3.2.
With each identified $\lambda_i$ we can return a list of extra informational items:
\begin{itemize}
  \item the number of elements validating ${_u\lambda}_i$ in the ${_uL}$ cluster,
  \item the number of elements validating ${^s\lambda_i}$ in the ${^sL}$ cluster,
  \item the actual radius of the ${_uL}$ cluster around ${_u\lambda}_i$,
  \item the actual radius of the ${^sL}$ cluster around ${^s\lambda}_i$,
\end{itemize}
The cardinality of the ${_uL}$ cluster ${_uC}_1$, which indicates how many
decimated analyses succeeded in retrieving $\lambda_i^u$, indicates the
level of validation of the retrieved $\lambda_i$, while that of the
${^sL}$ cluster ${^sC}_1$, in combination with its radius, reflects the correct or poor resolution from the aliasing.
The radius of ${_uC}_1$ on the one hand and of ${^sC}_1$ on the other, is a measure of the perturbation suffered by respectively $\lambda_i^u$ and $\lambda_i^s$.  
Clusters with few elements and large radii indicate that the conclusion
may be wrong because of the inherent noise. 
The total number of clusters detected in ${_uL}$ and validated in ${^sL}$, is automatically a good
estimate of the model order $n$, as pictured in Figure
\ref{cluster_explanation}.

\subsection{Outlier scenario}
It may happen that not all elements ${^s\lambda}_{i,k}$ tied to the
${_u\lambda}_{i,k}$ in a detected cluster ${_uC}_2$ belong to a cluster
${^sC}_2$. In that case the remote elements 
in ${^sL}$ are discarded and an estimate for ${^s\lambda}_i=\lambda_i^s$ is still the center of gravity of ${^sC}_2$. 
Here the number of decimated analyses validating $\lambda_i$ is different in ${_uL}$ and ${^sL}$.
\subsection{Collision scenario}
\begin{figure}
  \centering
\includegraphics[scale=0.5]{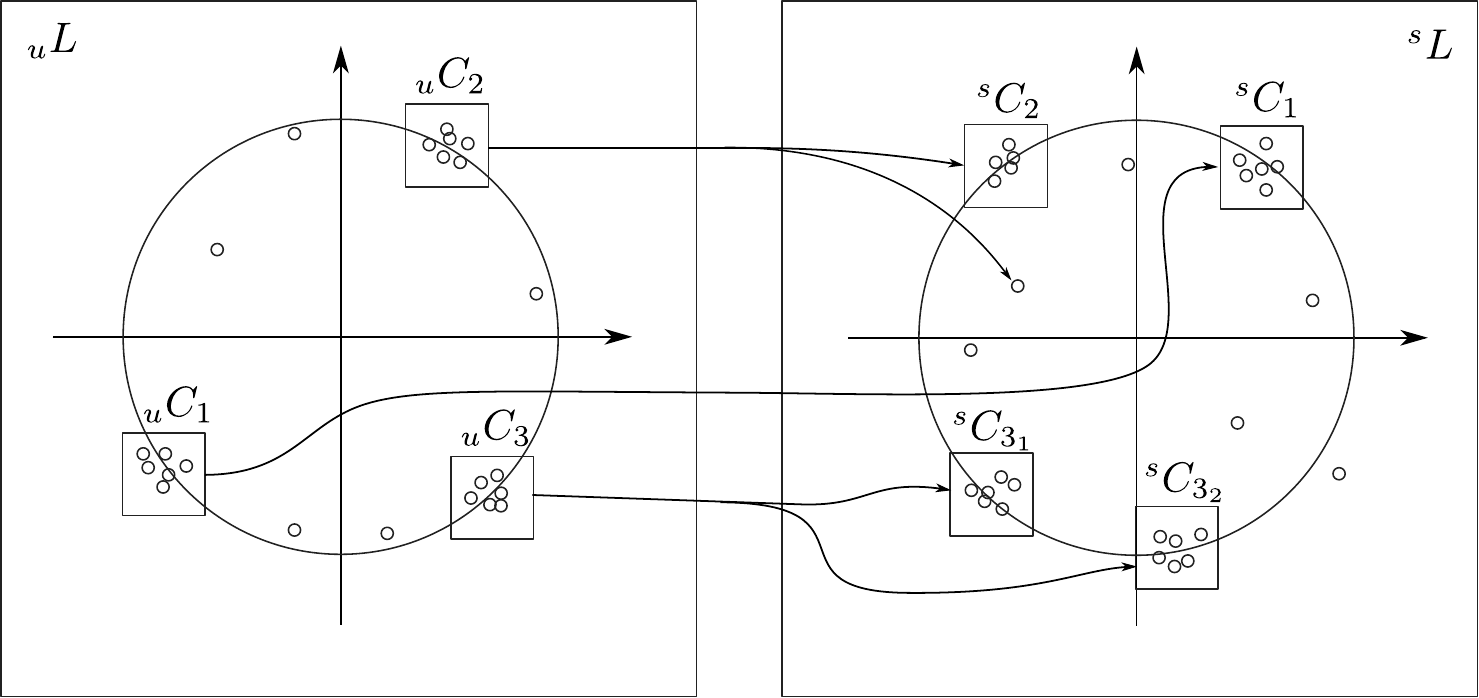}
\caption{\it The cluster algorithm as explained in Section 4, with the set ${_uL}$ at the left and the set ${^sL}$ at the right.} 
  \label{cluster_explanation}
\end{figure}

In cluster ${_uC}_3$ of Figure \ref{cluster_explanation}
a collision is involved. As in \eqref{msalpha_coll}
and pictured in Figure \ref{fig:collision}, the
${^{ms}\alpha}_{i,k}$ have identified more than one exponential
contribution. 
In ${^sL}$ different clusters of the ${^s\lambda}_{i,k}$ tied to the
${_u\lambda}_{i,k}$ in ${_uC}_3$ are identified instead of one large cluster. 
The centers of gravity of these individual clusters serve to identify the different generalized eigenvalues that have collided ${_uC}_3$ 
as a consequence of the aliasing.
\subsection{On the choice of $\delta$ and $m_\delta$}
At first sight, one may think that the new method, which now returns the
sparsity $n$ of $\phi(t)$ automatically by counting the number of detected
and validated clusters, heavily 
depends on a proper choice of the new parameters $\delta$ and $m_\delta$.  
This is actually not the case, because
usually multiple DBSCAN runs are performed, starting with a high validation 
rate $m_\delta \leq u$ and a small radius $\delta$, relaxing both 
gradually by decreasing $m_\delta$ and increasing
$\delta$, until the clusters detected in ${_uL}$ are not validated anymore 
by a cluster in ${^sL}$. Starting with more demanding validation
parameters identifies the most stable results first and then explores the
remaining results in a less strict way. We often choose $m_\delta$ in the
range $[0.75, 0.95] \times u$. 
As the ${_u\lambda}_i$ values are usually less affected by noise 
than the ${^s\lambda}_i$ which are obtained as solution of Vandermonde 
structured linear systems, the clusters in $_uL$ are generally denser than
the ones in $^sL$. So, in addition to the above, one should relax both
$m_\delta$ and $\delta$ a bit when moving with DBSCAN from $_uL$ to $^sL$.

\begin{figure}
  \centering
\includegraphics[width=5.3truecm]{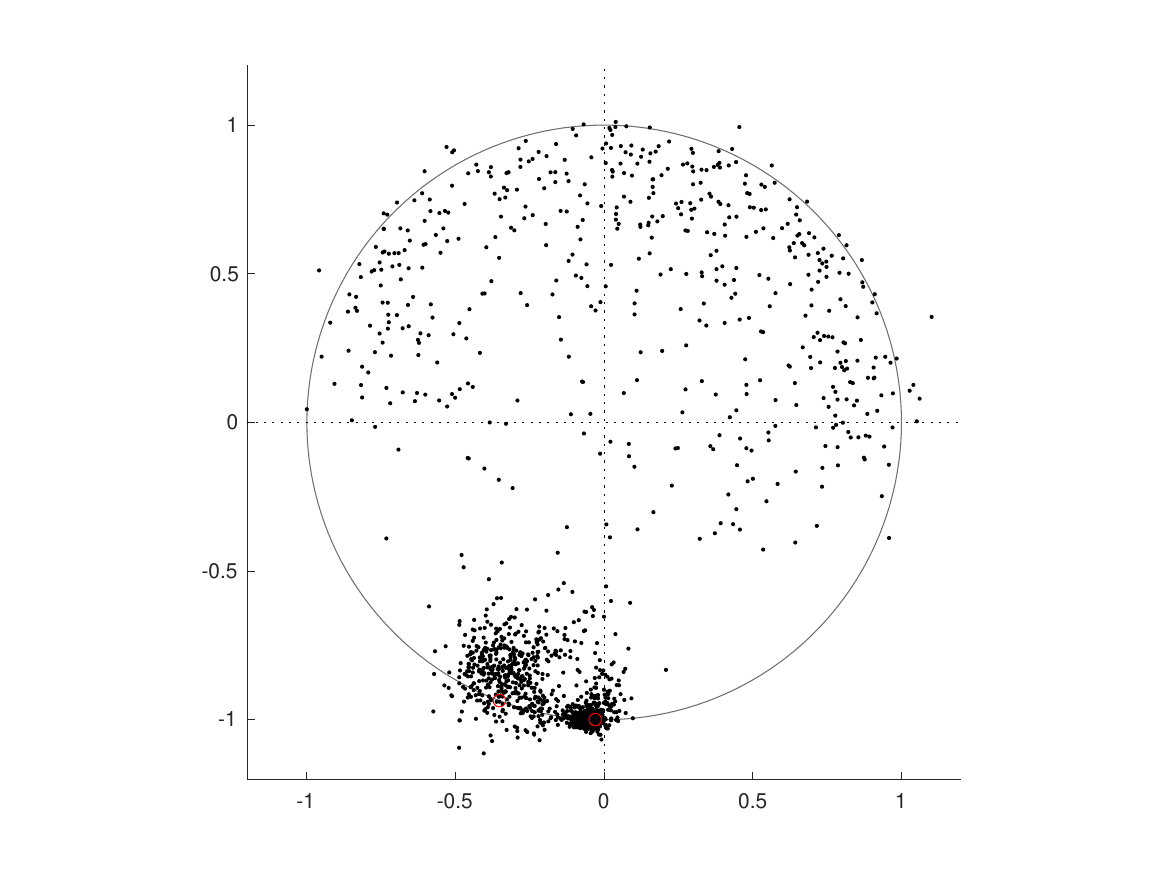}\quad
\includegraphics[width=5.3truecm]{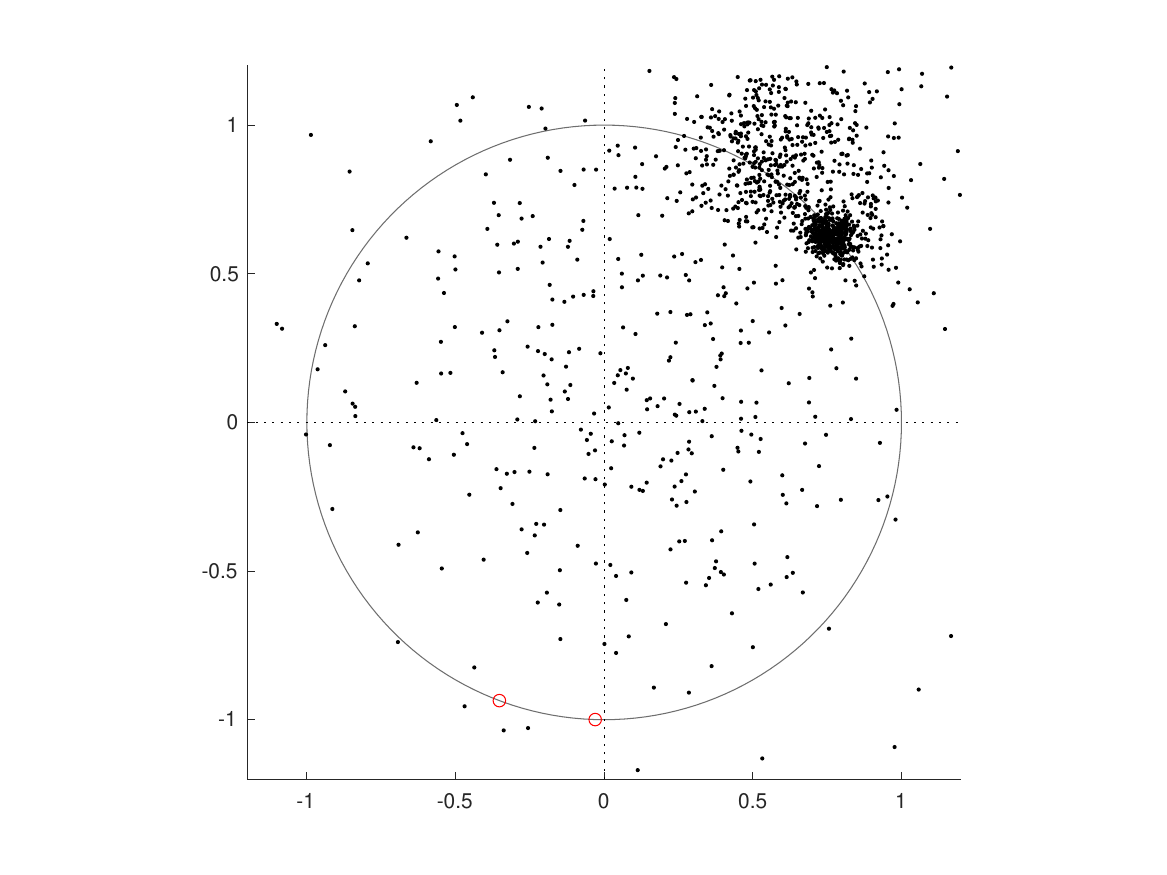} 
\caption{\it Joint $_7\lambda_{1,2}$ (at the left) and $^9\lambda_{1,2}$ 
(at the right) of \eqref{clustering} for all noise realizations.} 
\label{fig1-2}
\end{figure}

All the above is best illustrated with an example of an extreme case.
We take $n=2, \nu=3$ and $u=7, s=9$. 
We construct the illustration so as to generate difficult and large clusters, 
even exceeding the maximal theoretical cluster size $u$ 
of the decimation, by considering many different noise realizations. 
We aim for one dense and one diffuse cluster, with both getting
entangled as a consequence of their size and characteristic, and many unwanted
Froissart doublet poles on top. Consider for $j=0, \ldots, 13$ and
$\Delta=0.6028$, 
\begin{equation}
\phi(t_{s+ju})= 
\alpha_1 \exp(-\im 3.3576922 (s+ju)\Delta) + 
\alpha_2 \exp(\im 2.5206137 (s+ju)\Delta), 
\label{clustering}
\end{equation}
perturbed by white Gaussian noise with SNR $= 20$ dB.
Such a collection of $2 \times 14$ samples (for $u=7, s=0$ and $u=7, s=9$)
can result from the
described decimation technique. For the time being, let us call one such
collection a snaphot and let us
generate 512 snapshots by changing the noise realization in the
snapshot. 
From snapshot to snapshot 
the frequencies $\mu_i$ in \eqref{model}
remain unaltered. So the generalized eigenvalues $_u\lambda_i=
\exp(\mu_i (u\Delta)), i=1,2$
do not change from snapshot to snapshot. And neither do the values
$^s\lambda_i = \exp(\mu_i (s\Delta)), i=1,2$. 
In Figure \ref{fig1-2} 
one finds the $512 \times \nu = 1526$ values $_u\lambda_i$ at the left and
the same number of values $^s\lambda_i$ at the right. One clearly observes
two clusters in each of $_uL$ and $^sL$, 
as should be the case since $n=2$, one denser and one more
diffuse cluster, and a lot of scattered results coming from the fact that
$n$ is overestimated by $\nu>n$.

We now discuss the cluster detection in $_uL$ in more detail. Remember
that the process in $^sL$ is completely analogous, only with somewhat more
relaxed values for $m_\delta$ and $\delta$.

When running DBSCAN a single time, either with a small density $\delta$ or
a larger one, the correct result is not retrieved. When $\delta$ is small
(take $\delta=0.1$),
then only the denser cluster is revealed. When $\delta$ is larger, large
enough to go beyond the dense cluster (take $\delta=0.2$), 
then both clusters are joined into one. In both tries, we choose
$m_\delta=\lceil 0.85 \times 512 \rceil = 436$.

\begin{figure}
  \centering
\includegraphics[width=5.3truecm]{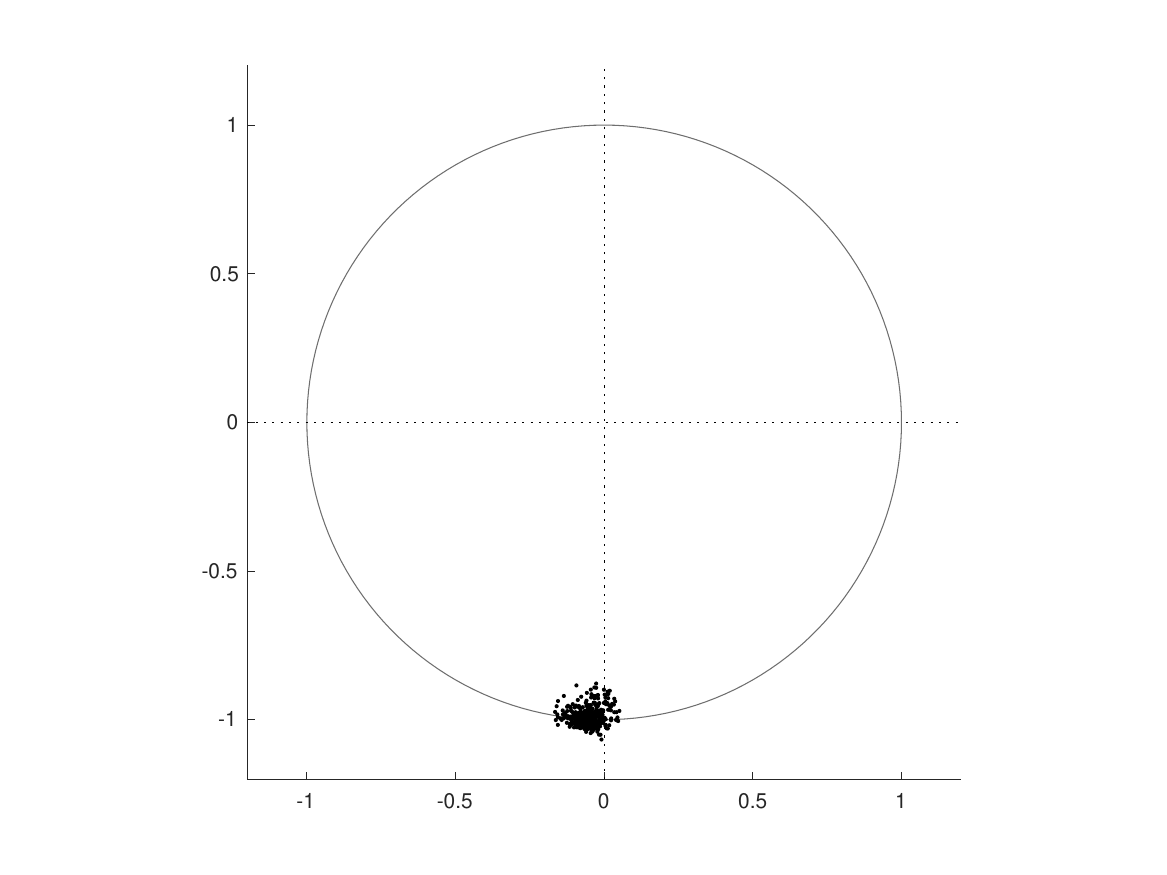}\quad
\includegraphics[width=5.3truecm]{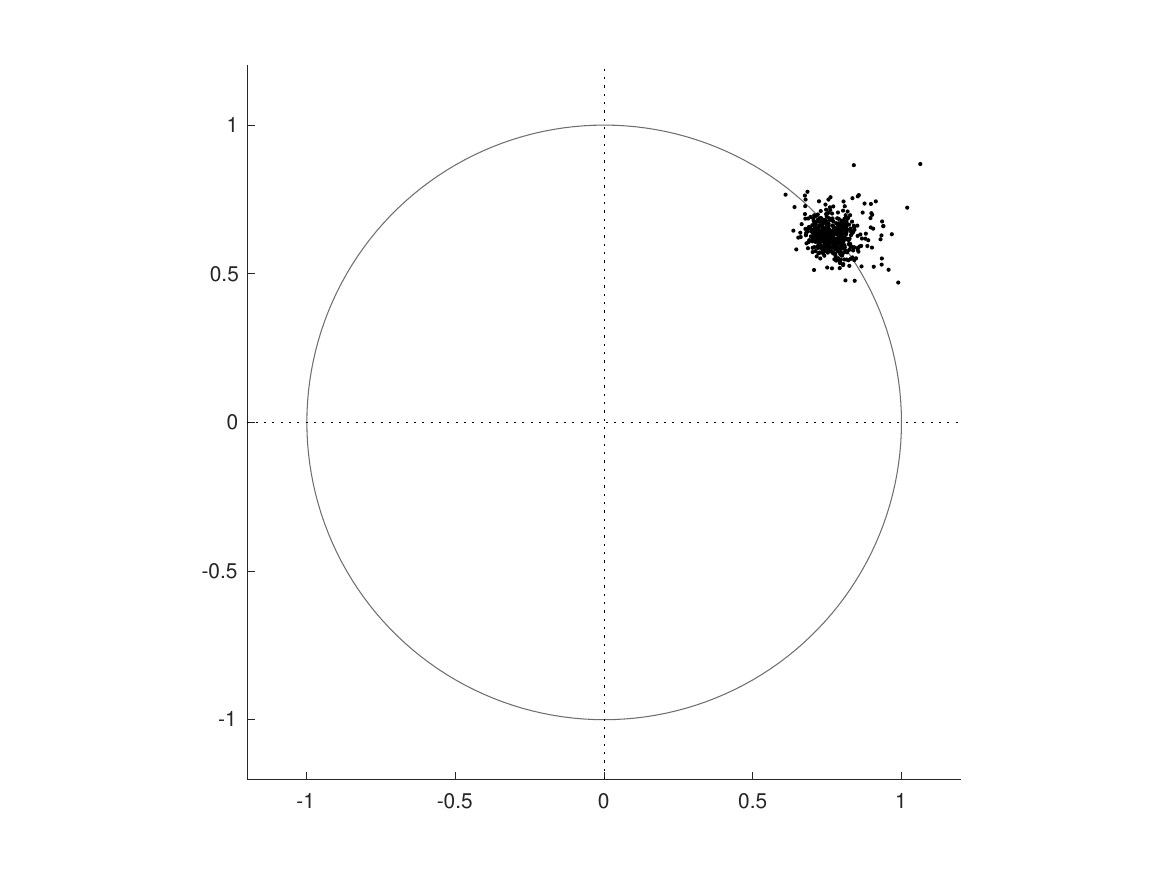} 
\caption{\it Dense cluster for \eqref{clustering}
retrieved in $_uL$ (left) and $^sL$ (right).}
\label{fig5-6}
\end{figure}

\begin{figure}
  \centering
\includegraphics[width=5.3truecm]{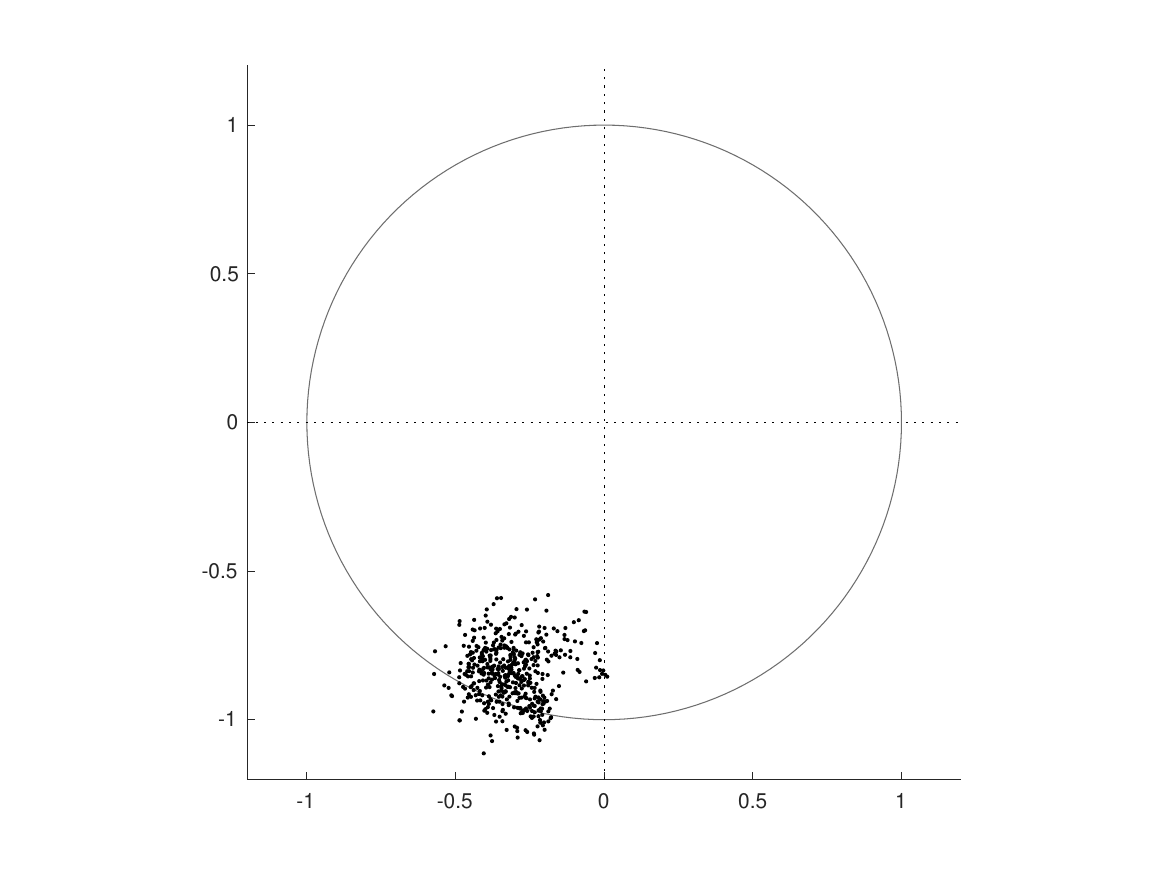}\quad
\includegraphics[width=5.3truecm]{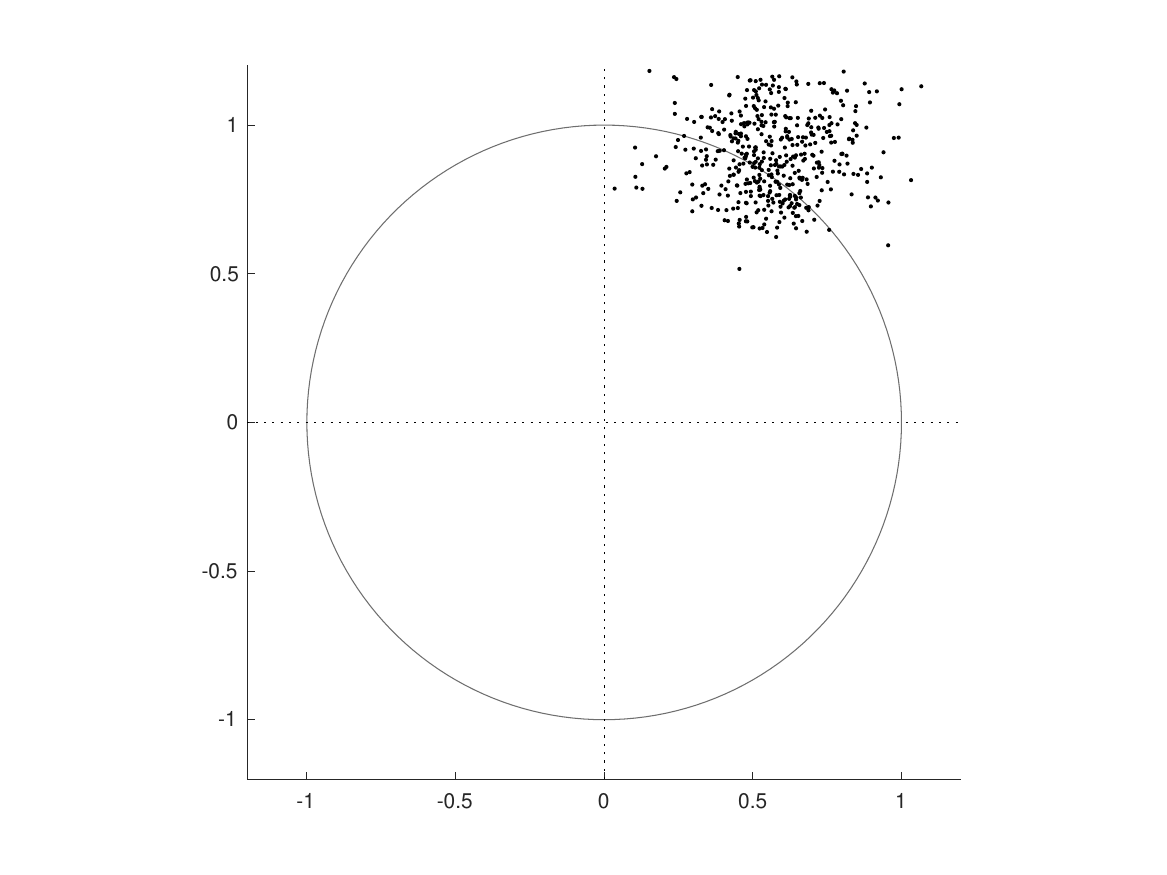} 
\caption{\it Diffuse cluster for \eqref{clustering}
retrieved in $_uL$ (left) and $^sL$ (right).}
\label{fig7-8}
\end{figure}

When using multiple DBSCAN runs, as explained above, the correct result is
revealed. Let us still fix $m_\delta=436$, but now vary $\delta$ from
small to large as $\delta=0.08 \ell, \ell=1, \ldots, 5$. 
With $\delta=0.08$ a first cluster is detected in respectively $_uL$ and
$^sL$, as shown in Figure \ref{fig5-6}. After detecting this
cluster, the concerned points are removed from the cluster analysis.
With $\delta=0.24$ a second cluster is detected in respectively $_uL$ and
$^sL$ and shown in Figure \ref{fig7-8}. For $\delta=0.16$ no cluster is
found. Actually, from $\delta=0.32$ on, up to $\delta=0.80$ no clusters 
are identified anymore. 

\begin{figure}
  \centering
\includegraphics[width=5.3truecm]{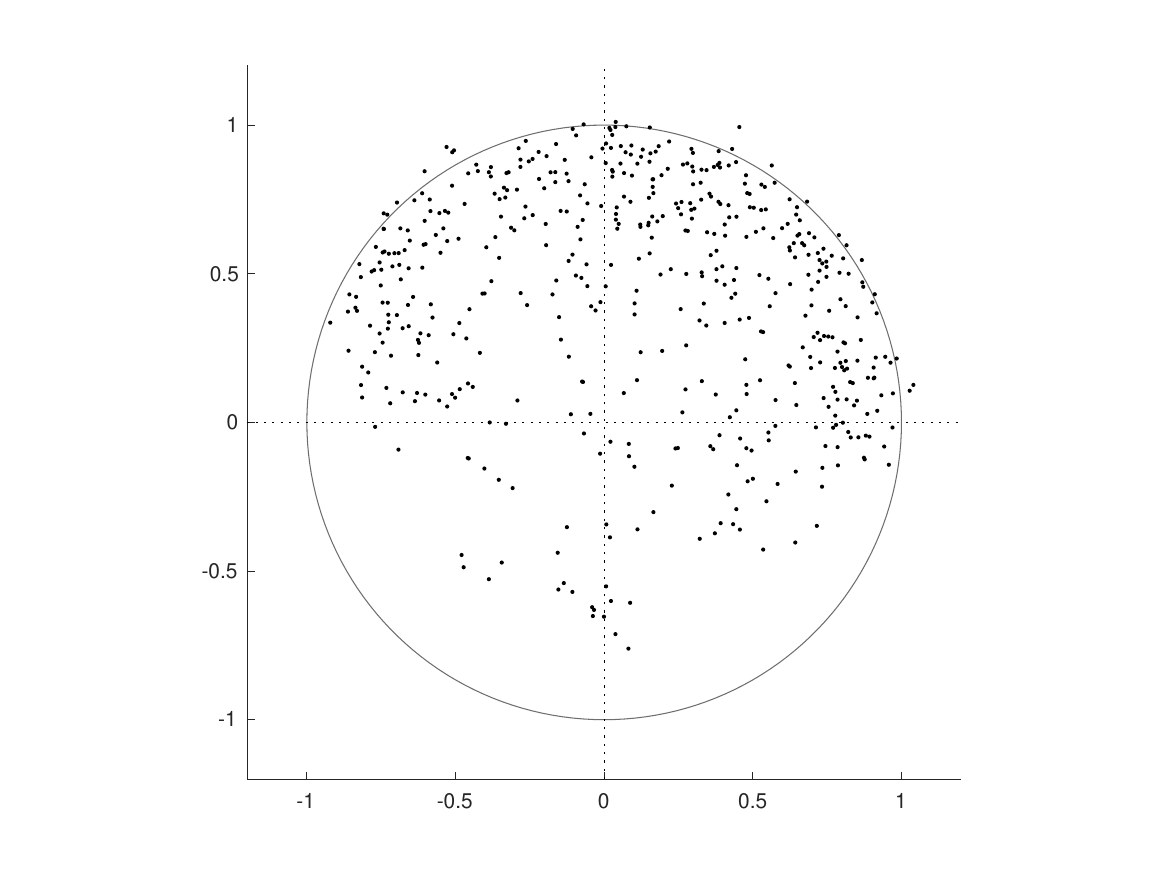}\quad
\includegraphics[width=5.3truecm]{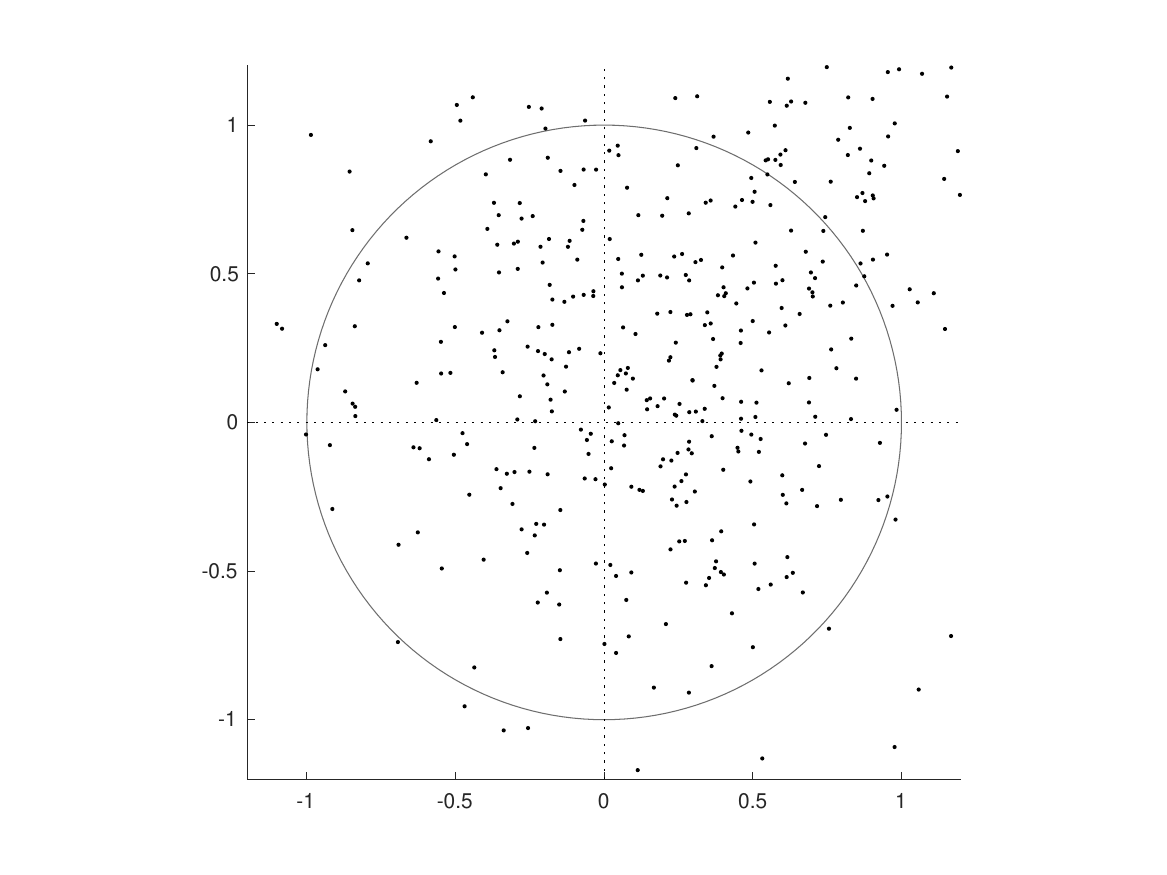} 
\caption{\it Candidate cluster for \eqref{clustering}
retrieved in $_uL$ (left) without associated cluster in $^sL$ (right).}
\label{fig9-10}
\end{figure}

When continuing the search with $\delta=0.88$ then a very diffuse cluster of 504
elements pops up again in $_uL$ (see Figure \ref{fig9-10} left) but
without confirmation by the clustering of the associated points in $^sL$
(see Figure \ref{fig9-10} right). So at this point, the search for clusters has
definitely been taken across reasonable values for $\delta$.

Another unrelated consideration that must be made with respect to $m_\delta$ is
the following. Let us denote $m_\delta = p \times u, 0 < p \le 1$, where
$p$ denotes the percentage of the maximal 
cluster size $u$ that we minimally
require for the cluster cardinality. If the set of $N$ signal samples still
contains some, say $\ell$, remaining outliers, 
despite an outlier filtering step which
is the initial preparation step prior to any further analysis of the data,
then in the worst case only $u-\ell$ elements can be found in each
cluster, and it is easily understood that we do not want $u-\ell$ to be
insignificantly small. 
This happens when each outlier contaminates a different subset 
$\Phi_k, k=1, \ldots, u$ of the decimation.
So we have the constraint $m_\delta \le u-\ell$ or $\ell \le (1-p)u$, which is
in practice too strict, but presents at least another indication of the
connection between $m_\delta$ and the quality of the data set. The strict
bound on $\ell$ actually guarantees that the decimation and clustering 
will not be bothered by the outliers.
When the strict bound is violated, the success of the method depends on the
location of the outliers and the probability that sufficient subsets
$\Phi_k, k=0, \ldots, u-1$ in the decimation are outlier free. 

  \section{Numerical illustration}
At this moment we introduce the acronym VEXPA for the new procedure that validates an exponential analysis carried out by a Prony-like method applied to each of the decimated signals.
In order to see the proposed method at work, we present the results of two
experiments, with the main aim to illustrate the extra features listed in
Section 4, which can now be added to whatever underlying
Prony-like method used for each separate decimated analysis.
For our experiments we use MP as the underlying method of choice to compute the aliased results ${_u\lambda}_{i,k}$ modelling the data $\Phi_k$ and to compute the recovery values ${^s\lambda}_{i,k}$ modelling the ${^{ms}\alpha}_{i,k}$. 
We then compare the VEXPA results to those of either the stand-alone MP
method or another popular Prony-like alternative \cite{St.Yi.ea:sta:94} and the atomic norm
minimisation (ANM) \cite{fastANM}.

All experiments are reproducible by downloading the matlab code and data
used in 5.1 and 5.2 
from {\tt cma.uantwerpen.be/publications}.

\subsection{Outlier experiment}
As discussed earlier, the cluster analysis makes the underlying exponential 
analysis algorithm more robust with respect to outliers that may have
escaped an outlier filtering step. In general, it is known that neither
methods of the Prony family nor basic implementations of ANM can deal
properly with impulsive noise or spikes in the data \cite{outliers}.

In Section 4.4 the relation $\ell \le u - m_\delta$,
between the undersampling factor $u$ of the
decimation, the number of outliers $\ell$ in the data and the choice for
the validation number $m_\delta$, is explained. We illustrate all this in
the following example. 

Consider $\phi(t)$ defined by the parameters $|\alpha_i|$,
$\arg(\alpha_i)$,
$\Im(\mu_i)$, $\Re(\mu_i)$, $i=1, 2, 3$ listed in Table 1.  
The total number of samples is $N=300$ and the bandwidth is $\Omega=1000$.  
White circular Gaussian noise with SNR $=30$ dB is added as well as some 
outliers.
We show the real part of the signal and disturb with real-valued outliers.
We notice no difference in the conclusions whether the outliers are real,
imaginary or complex. 
The following conclusions hold throughout.

\begin{figure}
  \centering
\includegraphics[width=6.0truecm]{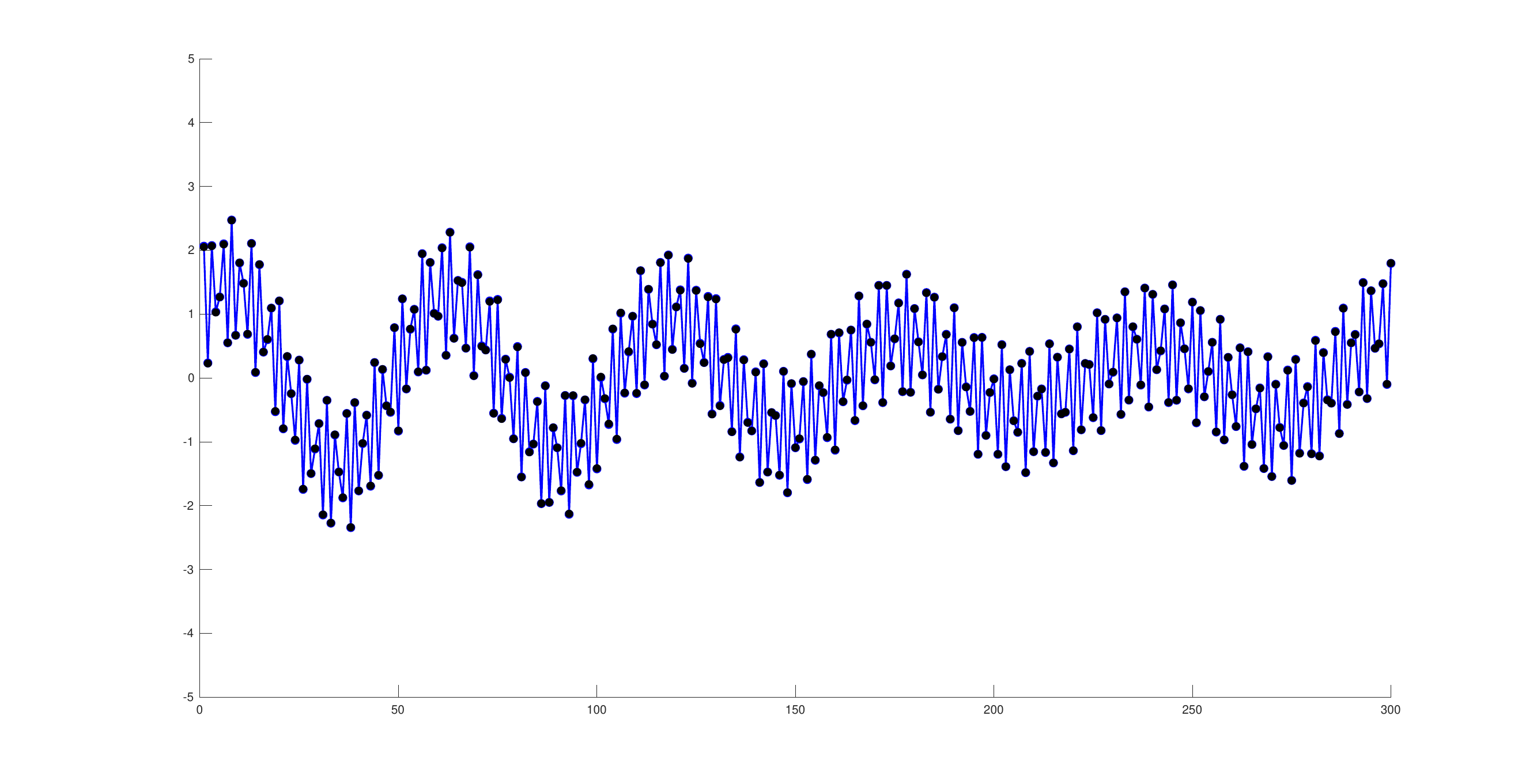}\quad
\includegraphics[width=6.0truecm]{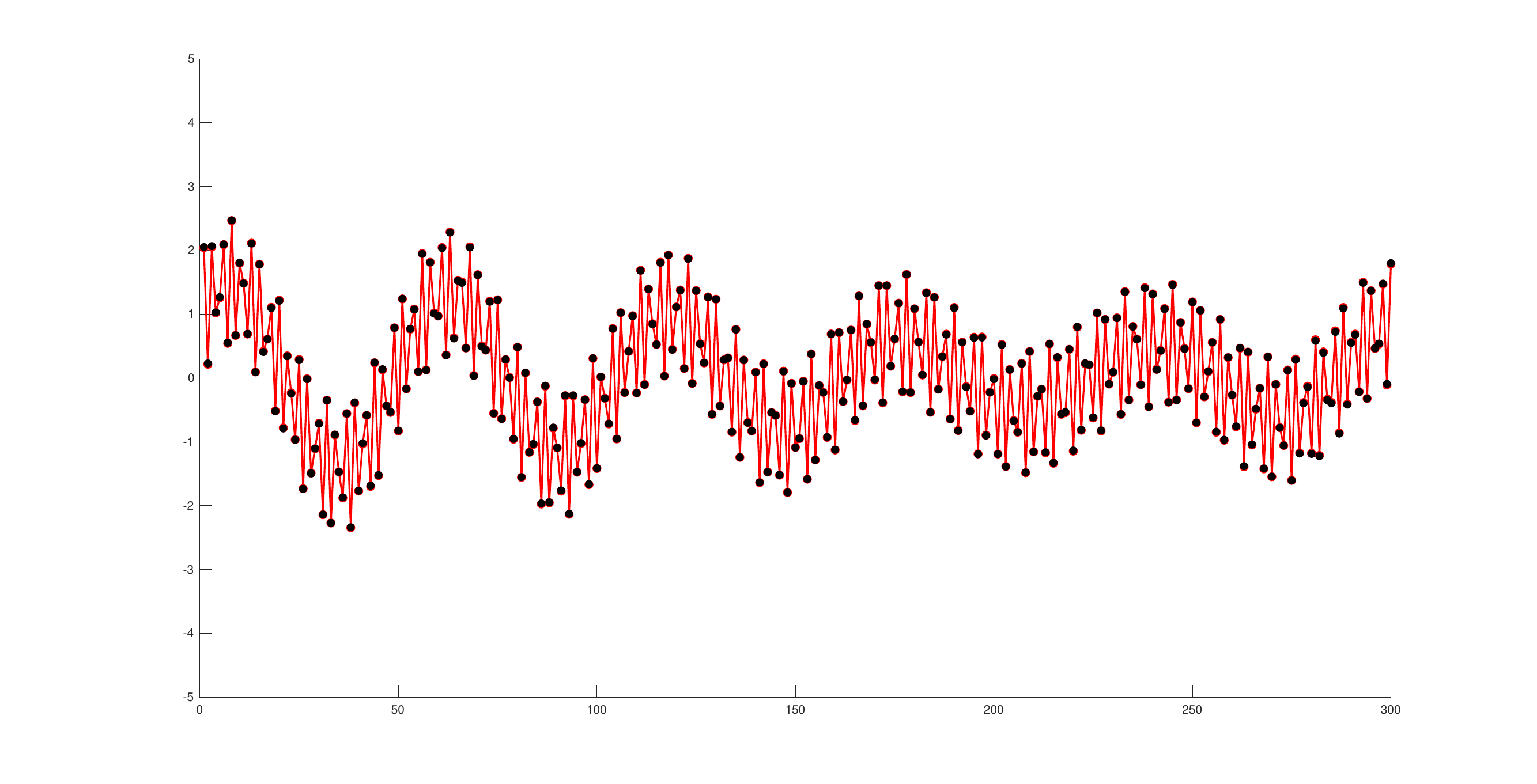} 
\caption{\it Outlier free reference with original data (in black), TLS-Prony
reconstruction (left, in blue) and VEXPA reconstruction (right, in red).}
\label{startoutlier}
\end{figure}

\begin{figure}
  \centering
\includegraphics[width=6.0truecm]{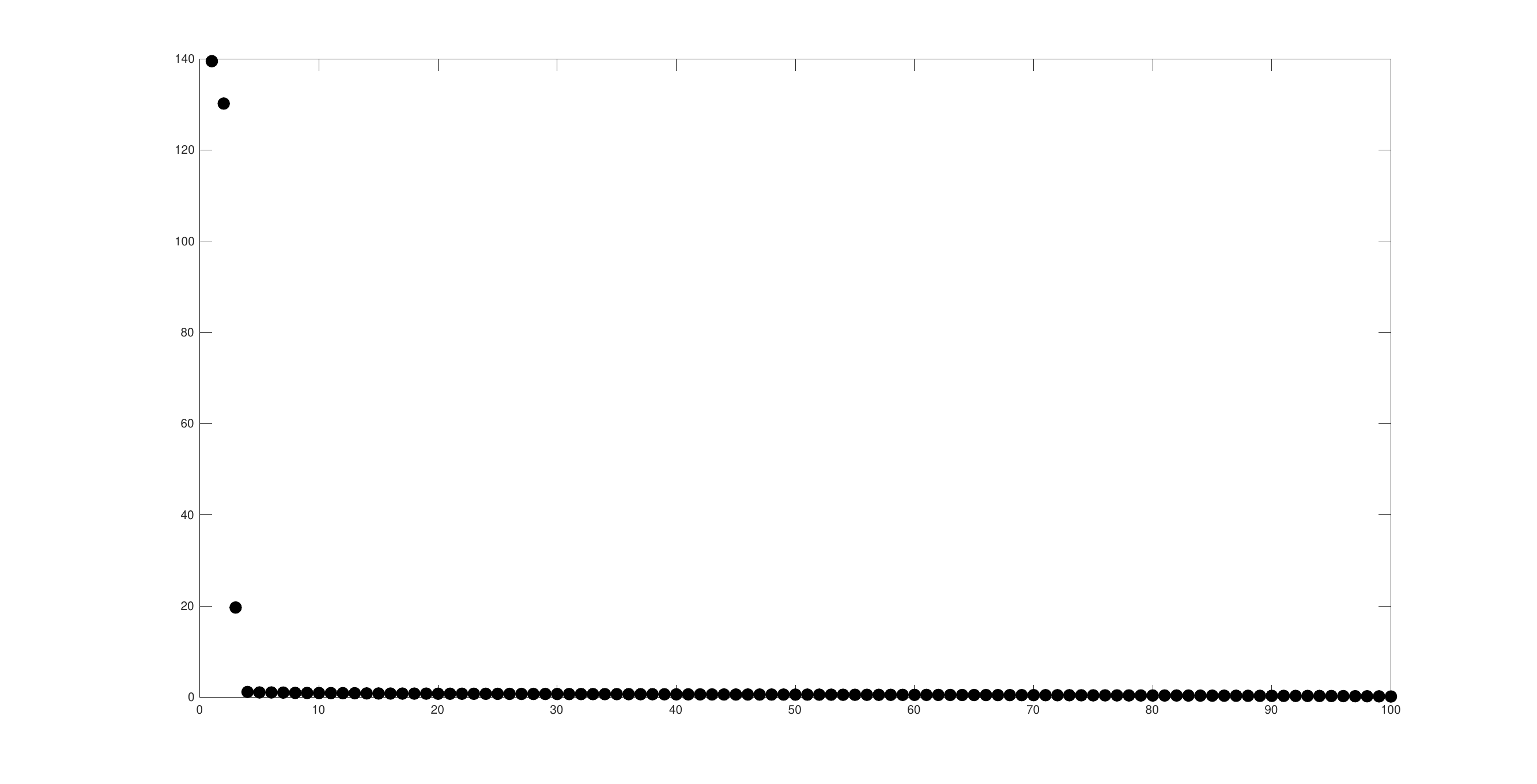}\quad
\includegraphics[width=6.0truecm]{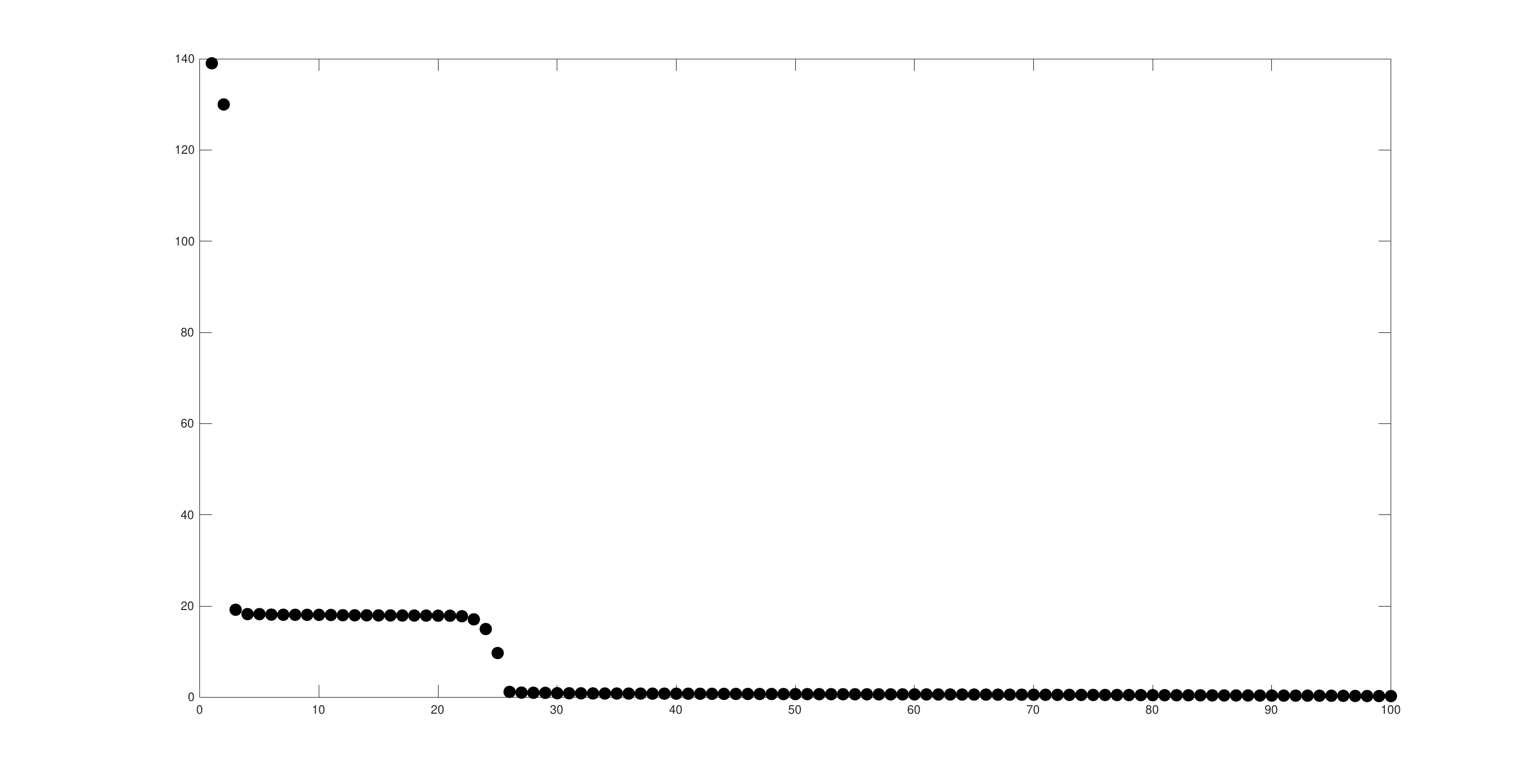} 
\caption{\it Numerical rank computation of $_1^0H_{200 \times 100}$, from
outlier free data (left) and single outlier data (right).} 
\label{SVDsec5}
\end{figure}

To establish some reference material, we first analyze the by noise corrupted 
but outlier free signal. 
On the one hand, we use a TLS-Prony method \cite{St.Yi.ea:sta:94} which
takes the numerical rank of the Hankel matrix $_1^0H_{(N-\nu) \times \nu}$ 
as a guess for the sparsity $n$. On the other hand, we compare this result 
to the output delivered by VEXPA with underlying the Prony-like algorithm MP
\cite{Hu.Sa:mat:90}.
Both methods recover the three terms and perform equally well. The
singular value
plot of the Hankel matrix $_1^0H_{200\times 100}$ used by the TLS-Prony
method, is given at the left in Figure \ref{SVDsec5}.
The root-mean-square errors are around 0.008 and the reconstructions are shown
in Figure \ref{startoutlier}, with the data in black, the TLS-Prony result
in blue at the left and the VEXPA result in red at the right.

\begin{figure}
  \centering
\includegraphics[width=6.0truecm]{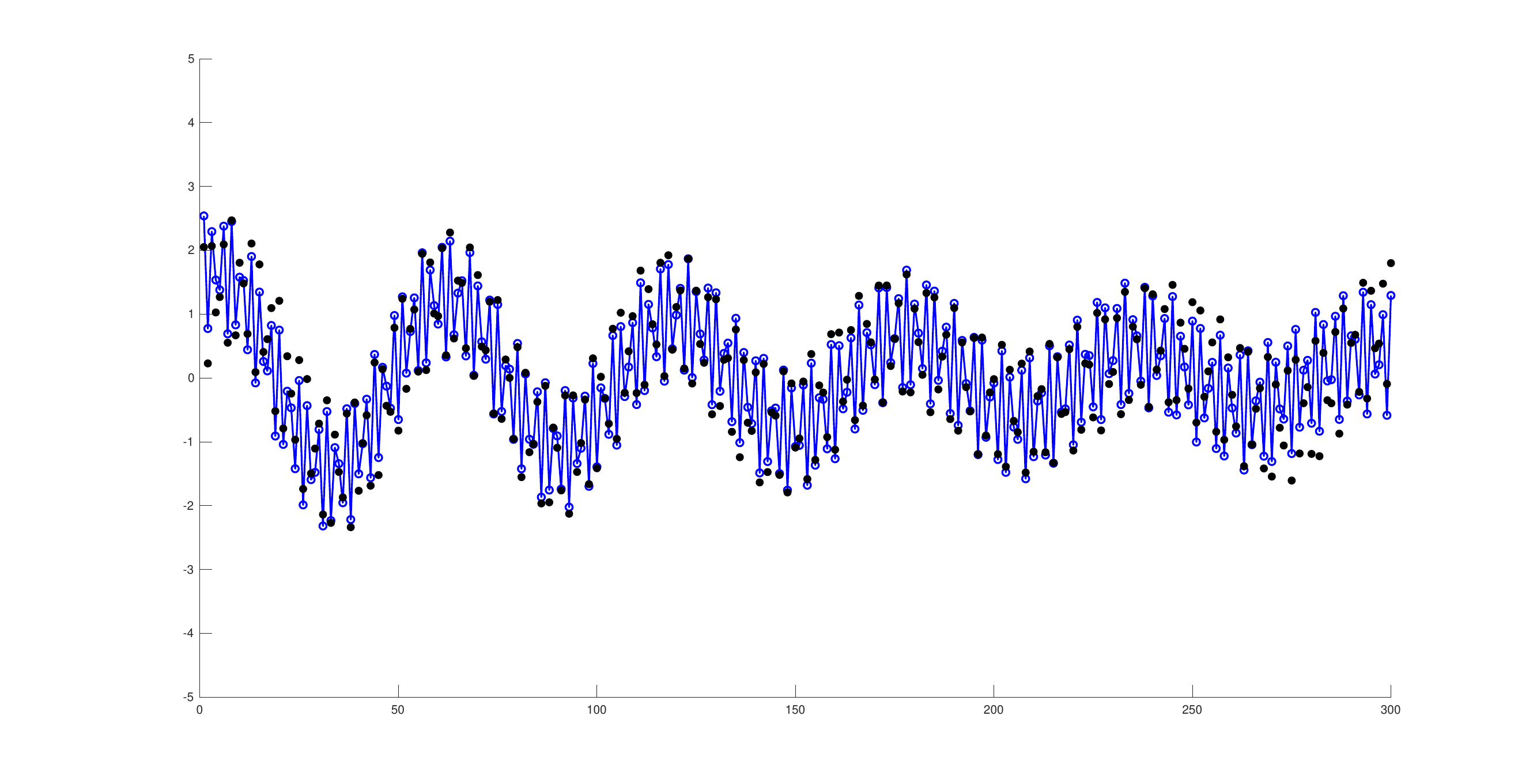}\quad
\includegraphics[width=6.0truecm]{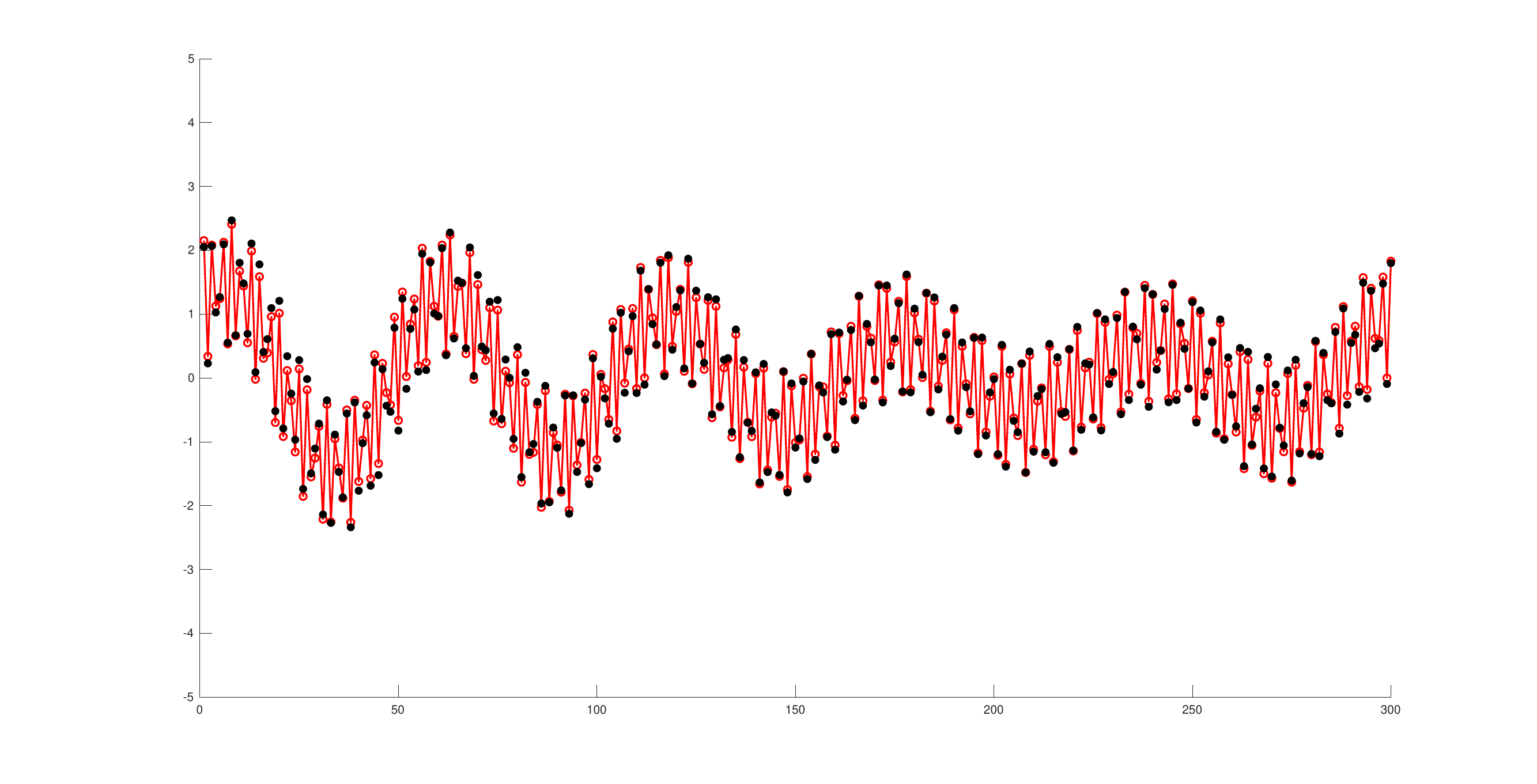} 
\caption{\it Outlier experiment with original data (in black), TLS-Prony
reconstruction (left, in blue) and VEXPA reconstruction (right, in red).}
\label{fig:outlierrecon}
\end{figure}

Already from one randomly placed outlier, we see that the RMSE of the
result delivered by VEXPA on top of an exponential analysis method is
generally less than that of the stand-alone method. 
Let us, for instance, subtract 18 from sample number 21. 
The singular value plot for the Hankel matrix
constructed with the data containing one outlier, is given in Figure
\ref{SVDsec5} at the right. Thresholding of the singular values of
${_1^0}H_{200\times 100}$ very obviously suggests to truncate all but two terms 
in the TLS step, although $n=3$ (opting for $n=25$ instead, results in
modelling the noise and outlier as well as the signal). 
So the third term is fully buried by the
outlier in the signal and we see this happen in all the subsequent
outlier tests on this signal. 
All subsequent singular value plots look similar and are
but slight variations of this one.
The reconstruction of the signal using the two recovered 
terms is shown in Figure
\ref{fig:outlierrecon}, at the left in blue. 
The reconstruction recovers quite well from the outliers,
but towards the end of the observation window, the signal
deviates more from the original. The RMSE, computed with
respect to the noisefree and outlier free signal, is 0.2912.

  \begin{table}
  \centering
    \label{table1}
  \begin{tabular}{|c|c|c|c|}
    \hline 
    $|\alpha_i|$ & $\arg(\alpha_i)$ & $\Im(\mu_i)$ & $\Re(\mu_i)$ \\
    \hline 
    1 & 0.3342 & $2\pi 417.764$ & -0.1 \\
    1 & 0.8084 & $-2\pi 17.4$ & 0 \\
    0.5 & 0.5880 & $-2\pi 19.5$ & 0 \\
    \hline
\end{tabular}
\caption{\it Section 5.1 experiment with $n=3$ and $N=300$.}
  \end{table}

\begin{figure}
    \centering
    \begin{minipage}{0.45\textwidth}
\includegraphics[width=\textwidth]{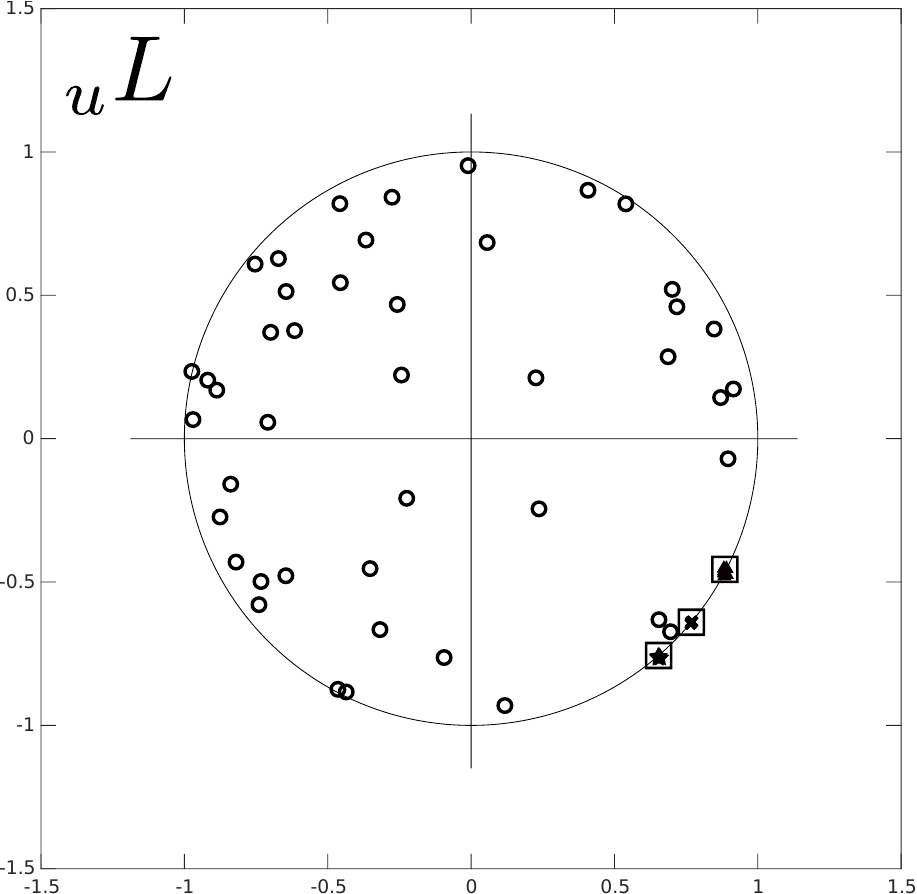}
    \end{minipage}
    \begin{minipage}{0.45\textwidth}
\includegraphics[width=\textwidth]{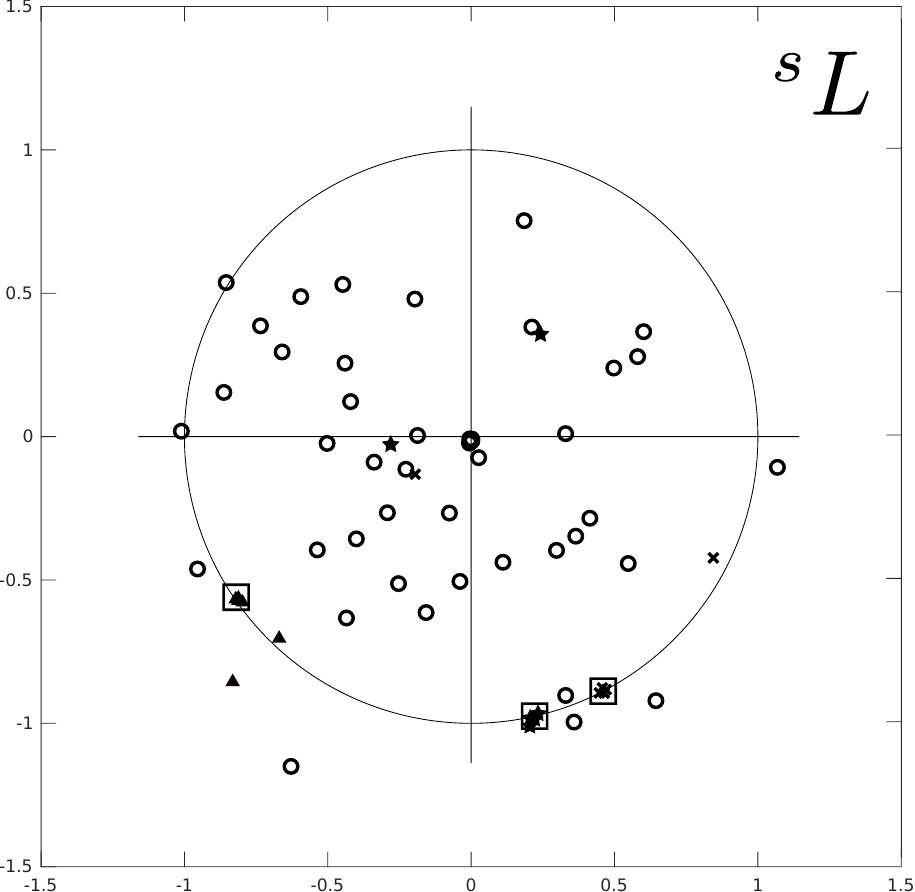}
    \end{minipage}
\caption{\it Cluster detection in ${_uL}$ (left) and ${^sL}$ (right) for the  outlier experiment.}
  \label{fig:outlierBu}
  \end{figure}

Next, the original signal is analyzed using VEXPA with underlying the
Prony-like algorithm MP without SVD thresholding, without a
guess for $n$. 
For VEXPA we take $u=7$ and $s=11$. 
So each $\Phi_k$ contains 42 or 41 samples. 
The decimation $\Phi_0$ contains the outlier.
So we can expect to find clusters of 6 elements in 
${_uL}$ instead of 7. Let us choose $m_\delta=5$.
We identify the ${_uL}$ clusters using increasing $\delta$-values, say
$\delta=0.01, 0.03, 0.05$, to isolate the most stable results first. 
In Figure \ref{fig:outlierBu} we show the results of the DBSCAN cluster 
analysis on ${_uL}$
and ${^sL}$. 
The VEXPA add-on clearly identifies $n=3$ exponential terms in the
signal and reconstructs the signal
quite reliably over the whole time interval. 
The signal reconstructed from the VEXPA output is depicted in Figure
\ref{fig:outlierrecon}, at the right in red. 
The RMSE is now 0.1164.

Very similar results are found with somewhat more outliers.
The bound $\ell \le u-m_\delta$ for the worst case, where each
outlier disturbs a different subset $\Phi_k$, is known to be too strict.
But remember that we assume to be dealing with only a few
outliers that may have escaped a filtering step. 
We now illustrate the capabilities of the new method in three different
experiments, all starting from the same signal and noise level:

\begin{itemize}
\item some statistical information on the RMSE in case of $\ell=2$ outliers, 
which fits the constraint $\ell\le u-m_\delta$ when $u=7, m_\delta=5$,
\item a more specific case belonging to these statistical data, 
where two outliers are located in each other's vicinity,
\item a typical situation where $\ell=5$ outliers are randomly placed, 
in this case fortunately without affecting all $\Phi_k$. 
\end{itemize}

\begin{figure}
  \centering
\includegraphics[height=5.5truecm]{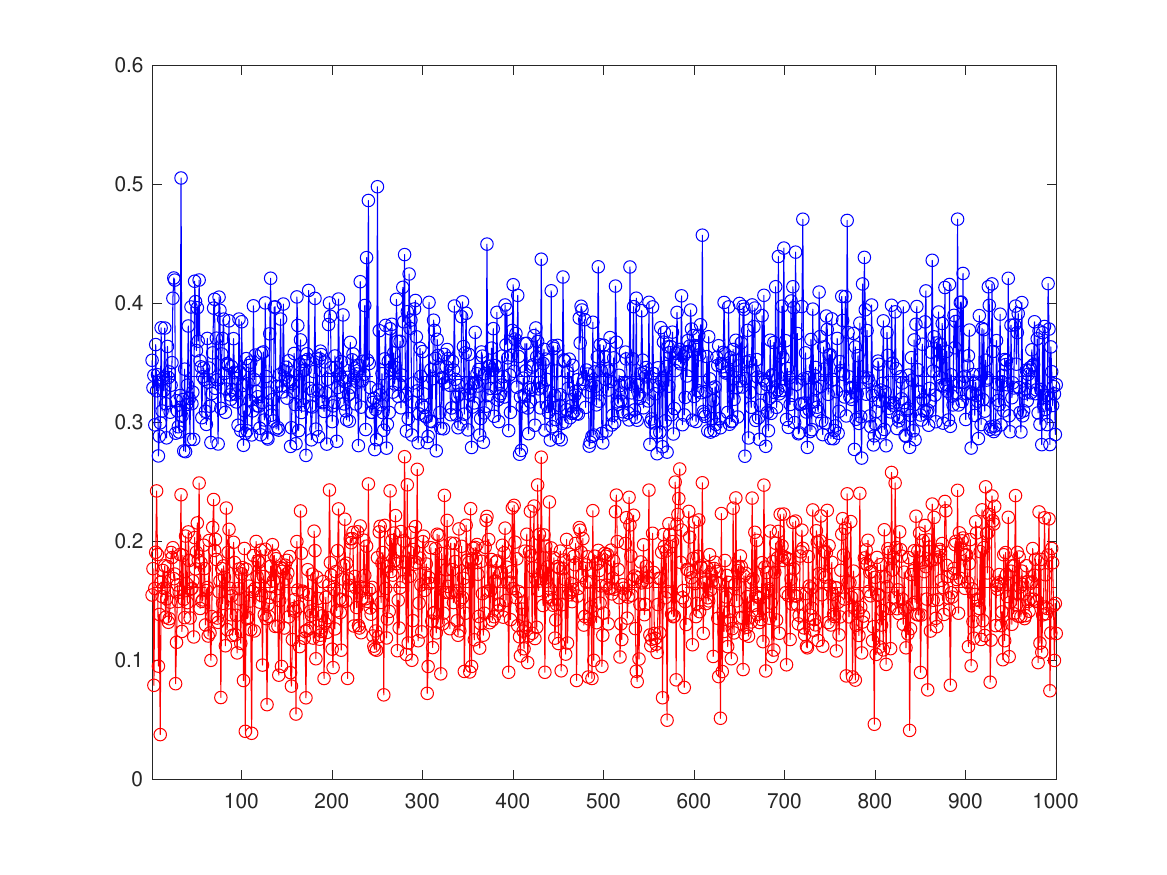}
\caption{\it RMSE of 1000 runs of TLS-Prony (blue, top) versus VEXPA (red, 
bottom) with 2 outliers randomly placed among the 300 samples.}
\label{RMSE}
\end{figure}

For the statistics, we add randomly chosen outliers belonging to the set 
$[-25,-15] \cup [15,25]$ to two randomly selected sample numbers
between 0 and 299. We do not change the noise, so that the effect that we
observe on the computation, compared to Figure \ref{startoutlier} with
reference RMSE values around 0.008, solely comes from the outliers. 
The experiment is repeated 1000 times. In Figure \ref{RMSE} we plot the
1000 RMSE of the TLS-Prony result versus that of the VEXPA result.

\begin{figure}
  \centering
\includegraphics[width=6.0truecm]{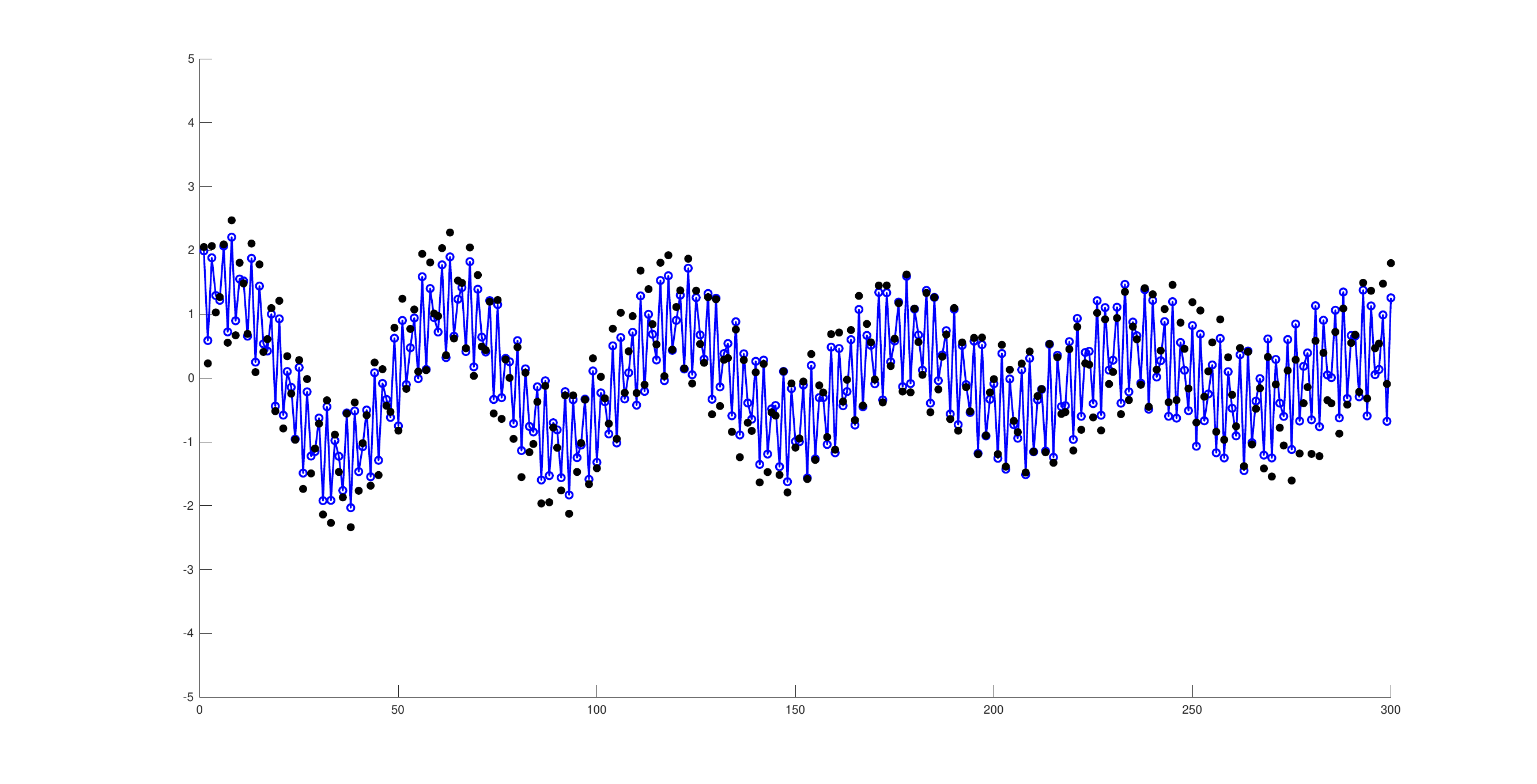}\quad
\includegraphics[width=6.0truecm]{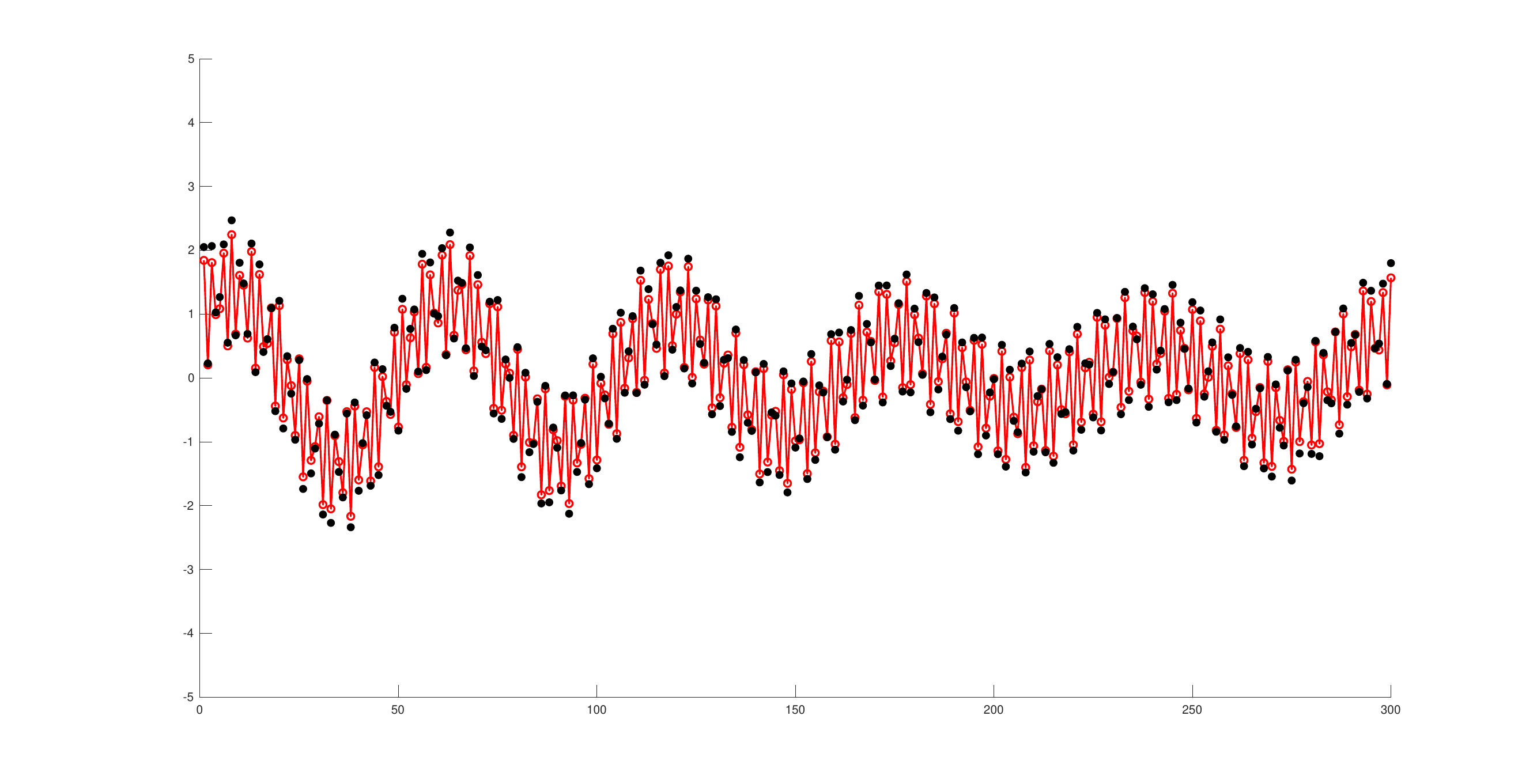} 
\caption{\it Outlier experiment with 2 outliers placed close to one
another.}
\label{2close}
\end{figure}

To illustrate the effect of 2 nearby outliers, we add respectively $-18$
and 24 to the sample numbers 21 and 25. Again the result is very similar. 
The TLS-Prony reconstruction with 2 terms is shown 
at the left in blue (see Figure \ref{2close}) with RMSE $=0.3241$.
The VEXPA reconstruction using 3 terms is shown at the right in red (see Figure
\ref{2close}) with RMSE $=0.1393$. 

\begin{figure}
  \centering
\includegraphics[width=6.0truecm]{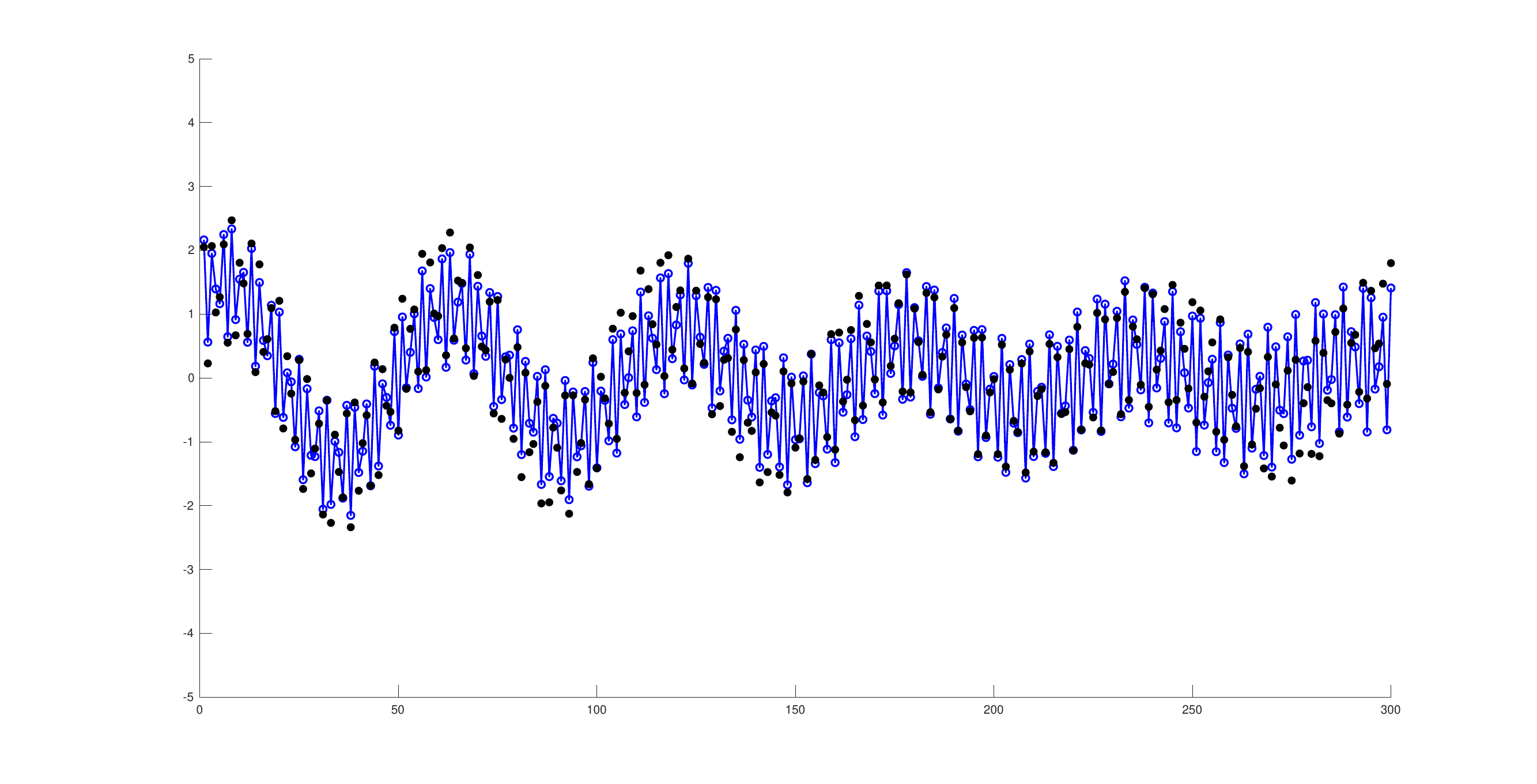}\quad
\includegraphics[width=6.0truecm]{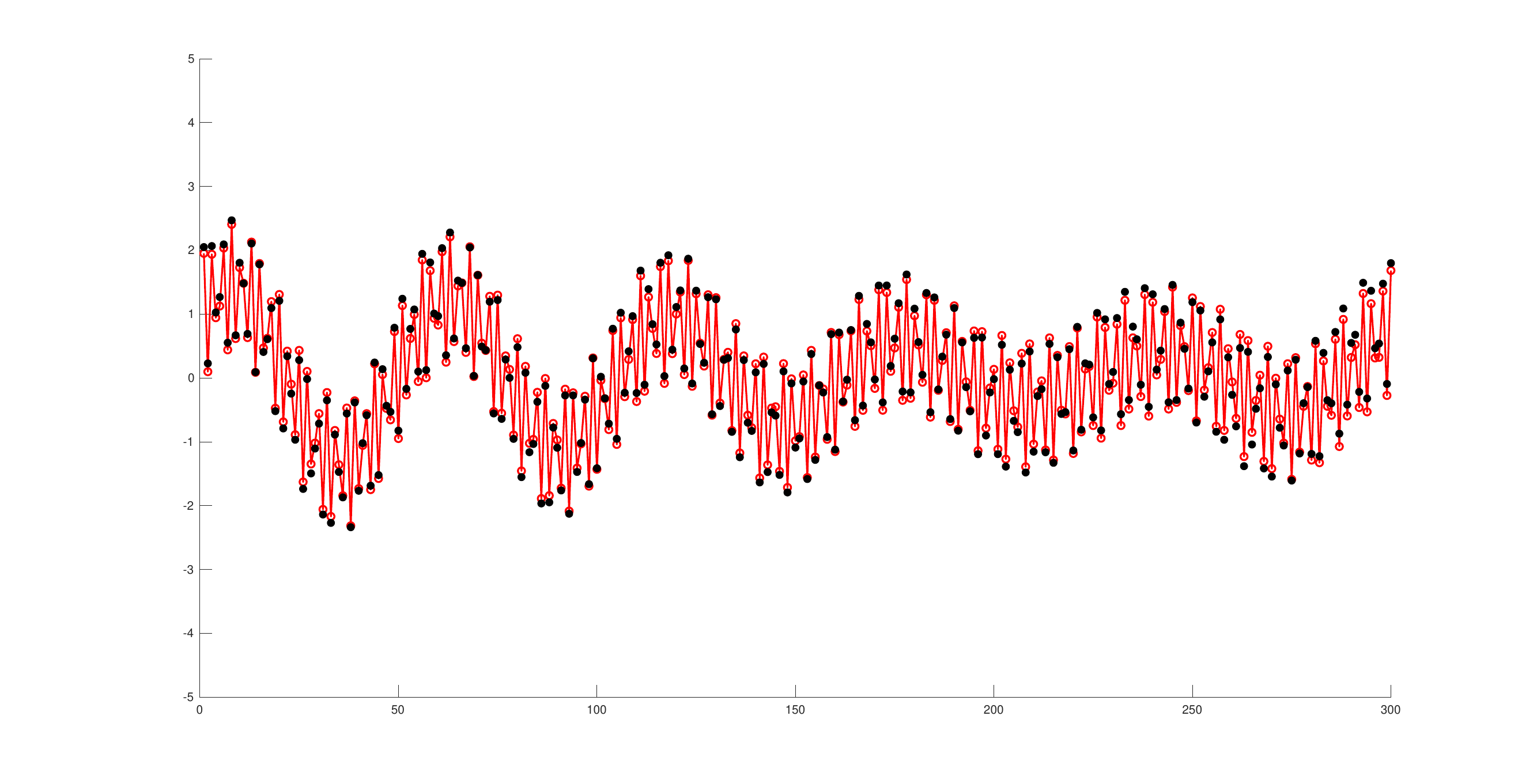} 
\caption{\it Outlier experiment with 5 outliers distributed over the
signal.}
\label{5outlier}
\end{figure}

To create 5 outliers we, for instance, 
respectively add $-18, 24, 17, -13, 20$ to the
sample numbers 21, 25, 134, 188, 258. The TLS-Prony method again suggests
very obviously to use two terms for the reconstruction (see Figure
\ref{5outlier} at the left in blue) with a RMSE $=0.3320$. The VEXPA method with the same
parameters for $\delta$ and $m_\delta$ reconstructs three terms (see
Figure \ref{5outlier} at the right in red) with a RMSE $=0.1390$.

How has the computation of the $\alpha_i$ profited from the cluster analysis as well? 
Since the clusters in ${_uL}$ consist of 5, 6 or 7 elements, we can deduce
precisely which subset(s) $\Phi_k$ 
did not contribute to the validation and so we
can omit all data points from such subset(s) in the linear system
delivering the parameters $\alpha_i$. 
So the computation of the $\alpha_i$ starts from outlier filtered data.  

\subsection{High noise experiment}
For our second experiment we consider a signal $\phi(t)$ defined by the
parameters $|\alpha_i|, \arg(\alpha_i), \Im(\mu_i), \Re(\mu_i), i=1, \ldots, 12$ in Table 2.
The total number of samples is again $N=300$, but now with $\Omega=100$.
We perturb the samples with white circular Gaussian noise of increasing SNR. 
The perturbed signal is then analysed using MP (with the added information
that $n=12$) on the one hand and
VEXPA (on top of MP without the added information that $n=12$) on the other. 
For the latter we choose $u=7$ and $s=6$.
We pass the correct model order $n$ only to MP. 
The new VEXPA add-on detects it automatically (for ${_u}L$ and ${^s}L$ we
respectively take $m_\delta=6$ and 4, and we choose $\delta=0.1$ twice).
For each SNR this experiment is repeated 500 times. The exponential
analysis using MP is fed a $200 \times 100$ generalized eigenvalue problem
which is being reduced to $n=12$ columns after performing an SVD step.
On the other hand, each individual decimation solves a $27 \times 15$ (or
$26 \times 15$) generalized eigenvalue problem and afterwards a cluster
analysis is performed on the combined outputs of the $u=7$ decimations.

\begin{figure}
  \centering
\includegraphics[scale=0.22]{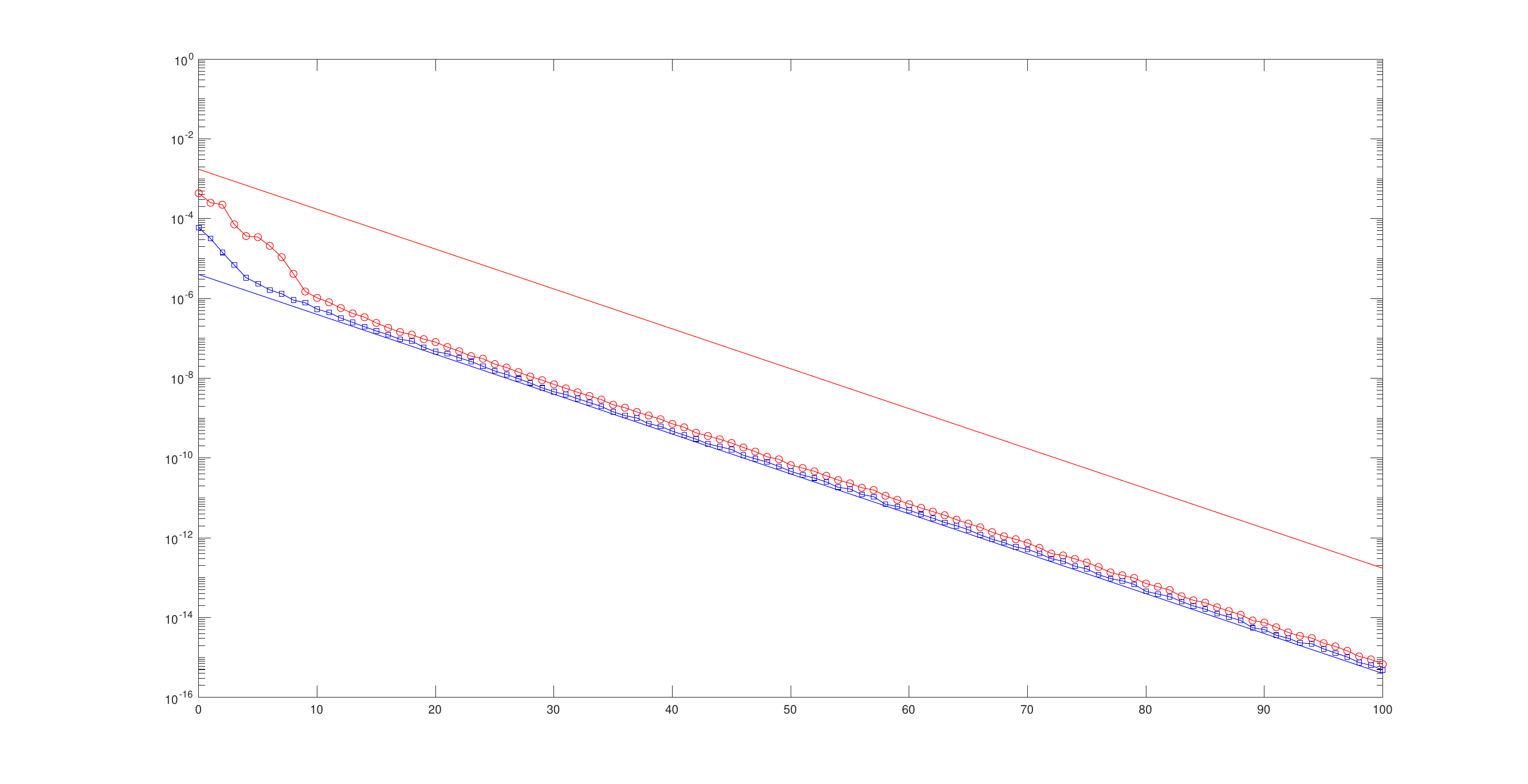}

\caption{\it Variance of MP (blue) and VEXPA (red), compared to the
Cramer-Rao lower bounds obtained as in Figure \ref{CRLBexample}.}
  \label{fig:CRLB}
\end{figure}

Up to SNR $=10$ dB both MP and VEXPA closely approach the desired CRLB,
as can be seen from Figure \ref{fig:CRLB}, where we show the CRLB for both
$\Omega=100, N=300$ (in blue) and $\Omega=100/7, N=42$ (in red), as in Figure
\ref{CRLBexample}. Remember that for each decimation
VEXPA is acually departing from the latter situation. But
after combining the different decimation results, and effectively also
using all the samples, the variance favourably compares to the CRLB. In
the meantime, several extra's have been picked up:
\begin{itemize}
\item As a consequence of the
decimation, the computational
complexity is greatly reduced because of the smaller independent
generalized eigenvalue problems. 
\item Hence the
numerical conditioning is improved and the analysis is parallellizable. 
\item In
addition, as already mentioned, the model order $n$ is an automatic
byproduct of the cluster analysis. 
\end{itemize}

\begin{figure}
\centerline{\includegraphics[scale=0.26]{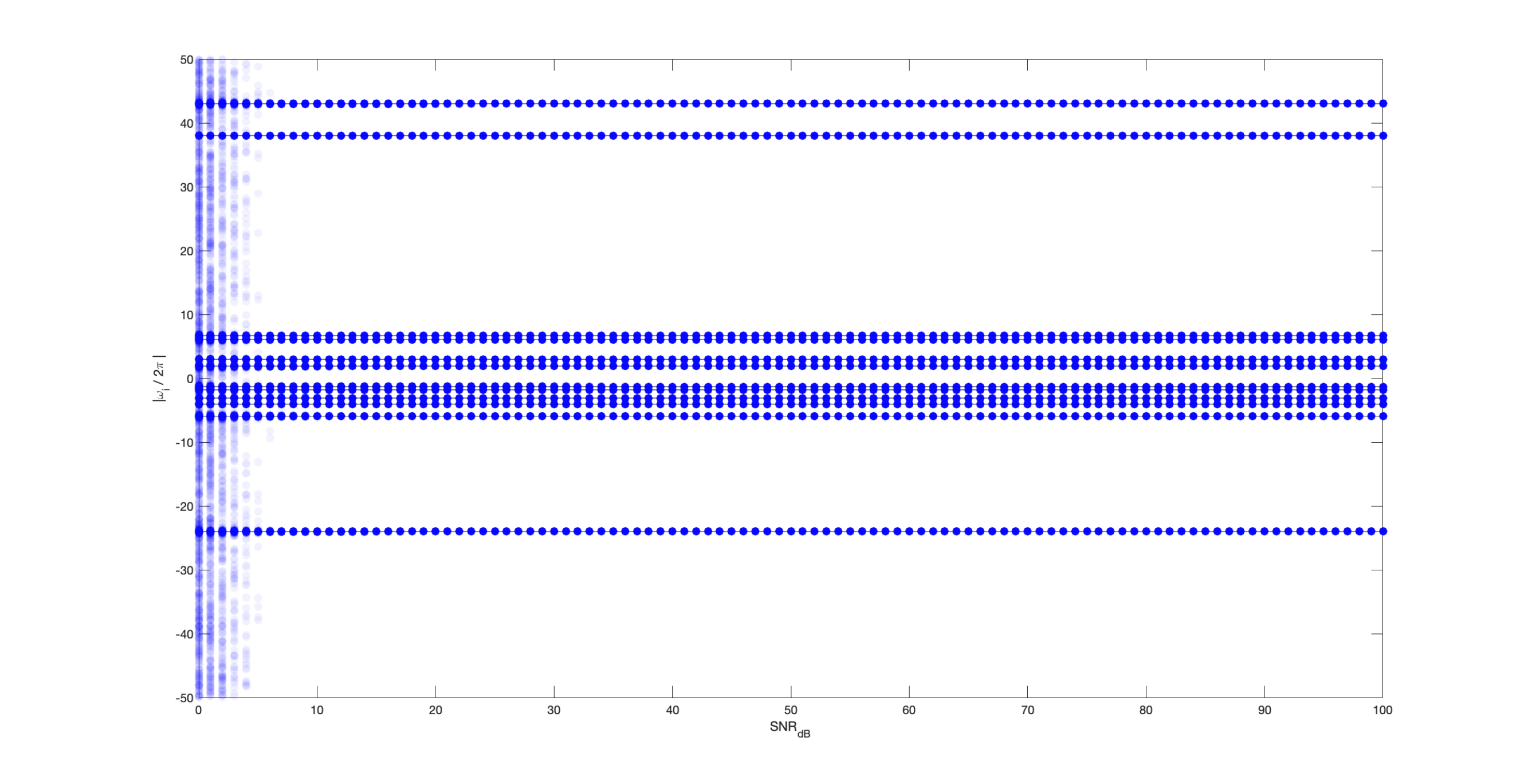}}
\centerline{\includegraphics[scale=0.26]{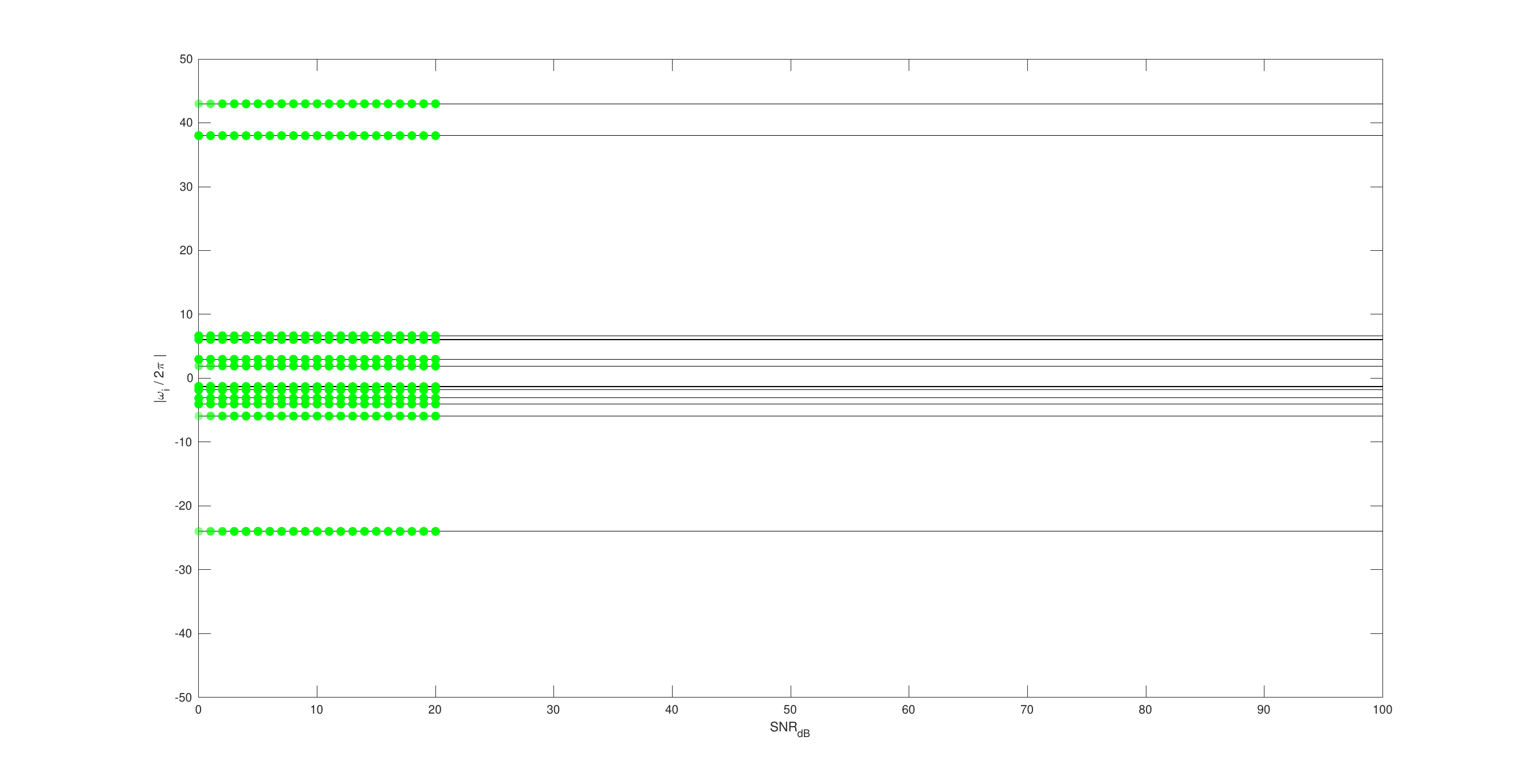}}
\centerline{\includegraphics[scale=0.26]{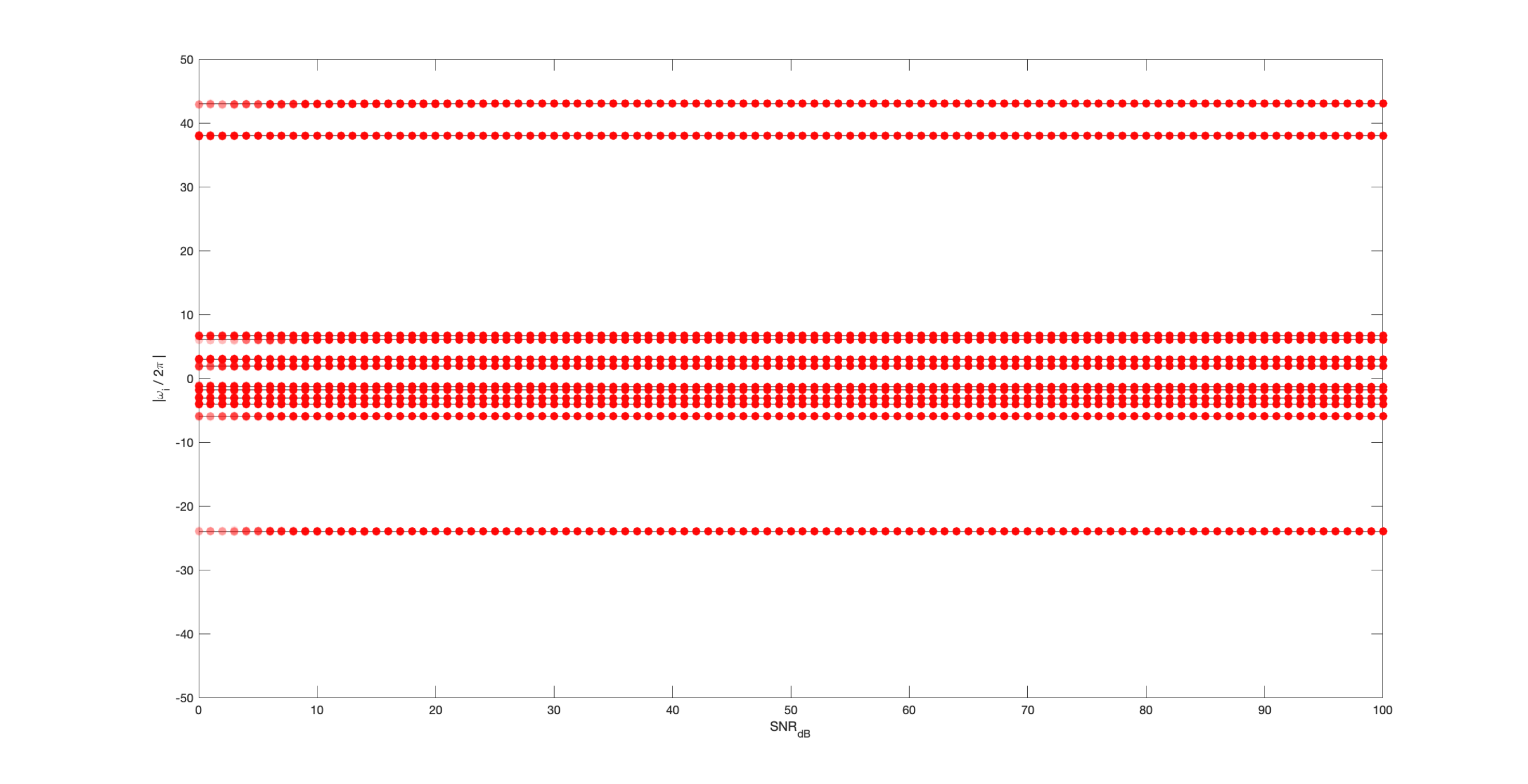}}
\caption{\it Retrieved $\Im(\mu_i)$ by MP (top, blue), ANM (middle, green)
and VEXPA (bottom, red).}
  \label{fig:poleslocations}
\end{figure}

Remains to discuss the gain in reliability. While Figure \ref{fig:CRLB}
on the comparison to the CRLB
is most interesting in the SNR interval $[10,100]$, the issue of
reliability becomes more fascinating for SNR values less than 10.
Since the signal is an undamped one, it can also be unravelled using the
ANM implementation in \cite{fastANM}. Despite the fact that this
implementation
enjoys an improved computational complexity of $O(N^2)$ per iteration step 
(maximum of 2000 iterations), it takes several hundred times longer per
execution than either MP or VEXPA. 
Therefore the method is only executed 100 times per SNR and
this for the more interesting interval of SNR values from 0 to 20.
The results of all runs, either 500 or 100, are superimposed in
Figure \ref{fig:poleslocations}: we show all retrieved $\Im(\mu_i)$-values for
MP (top), ANM (middle) and VEXPA (bottom).

For higher noise levels (smaller SNR) the stand-alone MP method
returns unreliable results, while the VEXPA method implemented on top of
MP detects when the signal is heavily perturbed, namely when 
fewer computed results are validated in the cluster analysis. 
So VEXPA, in its standard implementation, does not return
unreliable $\lambda_i$ output. When none of the results can be validated,
then VEXPA does not return $\lambda_i$ values at all. Such type of
reliability is also offered by ANM, as can be seen in the middle graph:
the retrieved frequencies are mostly correct, although some may be missing
in case of really small SNR. While the faster exponential analysis methods
of the Prony family traditionally suffer from an increased sensitivity to
noise (see the top graph on the stand-alone MP results), VEXPA adds as
good as the reliability that is normally offered by methods such as ANM
(see the middle graph of Figure \ref{fig:poleslocations}),
which can however be prohibitively slow. In addition, we point out that
VEXPA can easily be used on damped signals for which the ANM algorithm
does not qualify.

\begin{table}
\centering
\label{table2}
  \begin{tabular}{|c|c|c|c|}
    \hline
  $|\alpha_i|$ & $\arg(\alpha_i)$ & $\Im(\mu_i)$ & $\Re(\mu_i)$ \\
    \hline
  1 & 0 & $-2\pi 5.93$ & 0 \\
  2 & $\pi$ & $-2\pi 4.05$ & 0 \\
  2 & $\pi/4$ & $-2\pi 3.10$ & 0 \\
  2 & $\pi/8$ & $-2\pi 1.82$ & 0 \\
  2 & $3\pi/4$ & $-2\pi 1.31$ & 0 \\
  1 & $\pi/10$ & $2\pi 1.90$ & 0 \\
  3 & $-\pi$ & $2\pi 2.97$ & 0 \\
  1.5 & $-7\pi/8$ & $2\pi 6.05$ & 0 \\
  2 & 0 & $2\pi 6.67$ & 0 \\
  3 & $-78\pi/100$ & $2\pi 38$ & 0 \\
  1 & 0 & $2\pi 43$ & 0 \\
  1 & $\pi/5$ & $-2\pi 24$ & 0 \\
    \hline
\end{tabular}
\caption{\it Section 5.2 experiment with $n=12$ and $N=300$.}
\end{table}

\section{Conclusion}
Exponential analysis methods of the Prony type are more sensitive to noise. 
We offer an add-on technique that reconditions the problem statement and stabilizes and validates the computed results. 
As we illustrate in the numerical examples the algorithm works very well. 
In addition, the method
estimates the model order while performing the validation analysis.
The approach is highly suited for parallelization and hence further
improves the running time of the underlying Prony-like exponential
analysis, while offering a reliability comparable to that of the much more 
computationally intensive atomic norm minimization implementations.




\begin{thebibliography}{10}
\expandafter\ifx\csname url\endcsname\relax
  \def\url#1{\texttt{#1}}\fi
\expandafter\ifx\csname urlprefix\endcsname\relax\def\urlprefix{URL }\fi
\expandafter\ifx\csname href\endcsname\relax
  \def\href#1#2{#2} \def\path#1{#1}\fi

\bibitem{Ny:cer:28}
H.~Nyquist, Certain topics in telegraph transmission theory, Trans. Am. Inst.
  Electr. Eng. 47~(2) (1928) 617--644.
\newblock \href {https://doi.org/10.1109/T-AIEE.1928.5055024}
  {\path{doi: 10.1109/T-AIEE.1928.5055024}}.

\bibitem{Sh:com:49}
C.~E. Shannon, Communication in the presence of noise, Proc. IRE 37 (1949)
  10--21.


\bibitem{denoising}
B.~N. Bhaskar, G.~Tang, B.~Recht, Atomic norm denoising with applications to
  line spectral estimation, {IEEE} Transactions on Signal Processing 61~(23)
  (2013) 5987--5999.
\newblock \href {https://doi.org/10.1109/TSP.2013.2273443}
  {\path{doi: 10.1109/TSP.2013.2273443}}.

\bibitem{fastANM}
T.~L. Hansen, T.~L. Jensen, A fast interior-point method for atomic norm soft
  thresholding, Signal Processing 165 (2019) 7 -- 19.
\newblock \href {https://doi.org/10.1016/j.sigpro.2019.06.023}
  {\path{doi: 10.1016/j.sigpro.2019.06.023}}.

\bibitem{dePr:ess:95}
R.~de~Prony, Essai exp\'erimental et analytique sur les lois de la
  dilatabilit\'e des fluides \'elastiques et sur celles de la force expansive
  de la vapeur de l'eau et de la vapeur de l'alkool, \`a diff\'erentes
  temp\'eratures, J. Ec. Poly. 1 (1795) 24--76.

\bibitem{Ka.Le:ear:03}
E.~Kaltofen, W.-s. Lee, Early termination in sparse interpolation algorithms,
  J. Symbolic Comput. 36~(3-4) (2003) 365--400. 
\newblock \href
  {https://doi.org/10.1016/S0747-7171(03)00088-9}
  {\path{doi: 10.1016/S0747-7171(03)00088-9}}.

\bibitem{Ro.Ka:esp:89}
R.~Roy, T.~Kailath, {ESPRIT}-estimation of signal parameters via rotational
  invariance techniques, IEEE Trans. Acoust., Speech, Signal Process. 37~(7)
  (1989) 984--995.
\newblock \href {https://doi.org/10.1109/29.32276}
  {\path{doi: 10.1109/29.32276}}.

\bibitem{Hu.Sa:mat:90}
Y.~Hua, T.~K. Sarkar, Matrix pencil method for estimating parameters of
  exponentially damped/undamped sinusoids in noise, IEEE Trans. Acoust.,
  Speech, Signal Process. 38 (1990) 814--824.
\newblock \href {https://doi.org/10.1109/29.56027}
  {\path{doi: 10.1109/29.56027}}.

\bibitem{St.Yi.ea:sta:94}
W.~M. Steedly, C.-H.~J. Ying, R.~L. Moses,
  Statistical analysis of {TLS}-based {P}rony techniques, Automatica 30~(1) (1994) 
  115--129, special issue on statistical signal processing and control.
\newblock \href {https://doi.org/10.1016/0005-1098(94)90232-1}
  {\path{doi: 10.1016/0005-1098(94)90232-1}}.

\bibitem{Ch.Go:gen:06}
D.~Chu, G.~H. Golub, On a generalized eigenvalue problem for nonsquare pencils,
  SIAM J. Matrix Anal. Appl. 28~(3) (2006) 770--787.
\newblock \href {https://doi.org/10.1137/050628258}
  {\path{doi: 10.1137/050628258}}.

\bibitem{Cu.Le:how:18}
A.~Cuyt, W.-s. Lee,
  {How
  to get high resolution results from sparse and coarsely sampled data}, Appl.
  Comput. Harmon. Anal. 48~(3) (2020) 1066--1087.
\newblock \href {https://doi.org/10.1016/j.acha.2018.10.001}
  {\path{doi: 10.1016/j.acha.2018.10.001}}.

\bibitem{He:app:74}
P.~Henrici, Applied and computational complex analysis {I}, John Wiley \& Sons,
  New York, 1974.

\bibitem{Cu.Ts.ea:fai:18}
A.~Cuyt, M.~Tsai, M.~Verhoye, W.-s. Lee,
  {Faint
  and clustered components in exponential analysis}, Appl. Math. Comput. 327
  (2018) 93--103.

\bibitem{We.Mc:pro:63}
L.~Weiss, R.~McDonough, {P}rony's method, ${Z}$-transforms, and {P}ad\'e
  approximation, SIAM Rev. 5 (1963) 145--149.

\bibitem{Ba.My.ea:pad:89}
Z.~Bajzer, A.~C. Myers, S.~S. Sedarous, F.~G. Prendergast, Pad{\'e}-{L}aplace
  method for analysis of fluorescence intensity decay, Biophys. J. 56~(1)
  (1989) 79--93.

\bibitem{Nu:con:70}
J.~Nuttall, The convergence of {P}ad\'e approximants of meromorphic functions,
  J. Math. Anal. Appl. 31 (1970) 147--153.
\newblock \href
  {https://doi.org/10.1016/0022-247X(70)90126-5}
  {\path{doi: 10.1016/0022-247X(70)90126-5}}.

\bibitem{Po:pad:73}
C.~Pommerenke, {P}ad\'e approximants and convergence in capacity, J. Math.
  Anal. Appl. 41 (1973) 775--780.
\newblock \href
  {https://doi.org/10.1016/0022-247X(73)90248-5}
  {\path{doi: 10.1016/0022-247X(73)90248-5}}.

\bibitem{Ga:eff:72}
J.~Gammel, Effect of random errors (noise) in the terms of a power series on
  the convergence of the {P}ad\'e approximants, in: P.~Graves-Morris (Ed.),
  Pad\'e approximants, 1972, pp. 132--133.

\bibitem{Be:pad:96}
D.~Bessis, {P}ad\'e approximations in noise filtering, J. Comput. Appl. Math.
  66 (1996) 85--88.
\newblock \href
  {https://doi.org/10.1016/0377-0427(95)00177-8}
  {\path{doi: 10.1016/0377-0427(95)00177-8}}.

\bibitem{Gi.Pi:pad:97}
J.~Gilewicz, M.~Pindor, {P}ad\'e approximants and noise: a case of geometric
  series, J. Comput. Appl. Math. 87 (1997) 199--214.
\newblock \href
  {https://doi.org/10.1016/S0377-0427(97)00185-4}
  {\path{doi: 10.1016/S0377-0427(97)00185-4}}.

\bibitem{Gi.Pi:pad:99}
J.~Gilewicz, M.~Pindor, {P}ad\'e approximants and noise: rational functions, J.
  Comput. Appl. Math. 105 (1999) 285--297.
\newblock \href
  {https://dx.doi.org/10.1016/S0377-0427(99)00041-2}
  {\path{doi: 10.1016/S0377-0427(99)00041-2}}.

\bibitem{BARONE2005224}
P.~Barone,
  {On
  the distribution of poles of {P}ad{\'e} approximants to the {Z}-transform of
  complex {G}aussian white noise}, Journal of Approximation Theory 132~(2)
  (2005) 224 -- 240.
\newblock \href {https://dx.doi.org/10.1016/j.jat.2004.10.014}
  {\path{doi: 10.1016/j.jat.2004.10.014}}.

\bibitem{perotti2014identification}
L.~Perotti, T.~Regimbau, D.~Vrinceanu, D.~Bessis, Identification of
  gravitational-wave bursts in high noise using {P}ad{\'e} filtering, Physical
  Review D 90~(12) (2014) 124047.

\bibitem{Go.Mi.ea:sta:99}
G.~Golub, P.~Milanfar, J.~Varah, A stable numerical method for inverting shape
  from moments, SIAM J. Sci. Comput. 21 (1999) 1222--1243.

\bibitem{Be.Go.ea:num:07}
B.~Beckermann, G.~Golub, G.~Labahn, On the numerical condition of a generalized
  {H}ankel eigenvalue problem, Numer. Math. 106~(1) (2007) 41--68.
\newblock \href {https://dx.doi.org/10.1007/s00211-006-0054-x}
  {\path{doi: 10.1007/s00211-006-0054-x}}.

\bibitem{Kay:1993:FSS:151045}
S.~M. Kay, Fundamentals of Statistical Signal Processing: Estimation Theory,
  Prentice-Hall, Inc., Upper Saddle River, NJ, USA, 1993.

\bibitem{Yao1995CramerRaoLB}
Y.-X. Yao, S.~M. Pandit, {C}ramer-{R}ao lower bounds for a damped sinusoidal
  process, IEEE Trans. Signal Processing 43 (1995) 878--885.

\bibitem{YILMAZER2008561}
N.~Yilmazer, S.~Ari, T.~K. Sarkar,
  {Multiple
  snapshot direct data domain approach and {ESPRIT} method for direction of
  arrival estimation}, {D}igital {S}ignal {P}rocessing 18~(4) (2008) 561 --
  567.
\newblock \href {https://dx.doi.org/10.1016/j.dsp.2007.07.004}
  {\path{doi: 10.1016/j.dsp.2007.07.004}}.

\bibitem{Ester96adensity-based}
M.~Ester, H.-P. Kriegel, J.~Sander, X.~Xu, A density-based algorithm for
  discovering clusters in large spatial databases with noise, in: {KDD}'96
  Proceedings of the Second International Conference on Knowledge Discovery and
  Data Mining, {KDD}-96, AAAI Press, 1996, pp. 226--231.


\bibitem{outliers}
C.~Fernandez-Granda, G.~Tang, X.~Wang, L.~Zheng, {Demixing sines and spikes:
  Robust spectral super-resolution in the presence of outliers}, Information
  and Inference: A Journal of the IMA 7~(1) (2018) 105--168.
\newblock \href {https://dx.doi.org/10.1093/imaiai/iax005}
  {\path{doi: 10.1093/imaiai/iax005}}.


\end{thebibliography}





\end{document}